\documentstyle[amsfonts, amssymb, latexsym, epsfig, epic, eepic ,12pt]{amsart}

\newtheorem{theorem}{Theorem}[section]
\newtheorem{lemma}[theorem]{Lemma}
\newtheorem{proposition}[theorem]{Proposition}
\newtheorem{corollary}[theorem]{Corollary}
\newtheorem{definition}[theorem]{Definition}

\def\C{{\mbox{\rm\kern.24em
\vrule width.03em height1.43ex depth-.052ex \kern-.26em C}}}
\def\QSet{\mbox{\rm\kern.24em
\vrule width.03em height1.48ex depth-.051ex \kern-.26em Q}}
\def\Z{{\bf Z}}
\def\R{{\mbox{\rm I\kern-.22em R}}}
\def\P{{\bf P}}
\def\Q{{\bf Q}}
\def\T{{\bf T}}
\def\D{{\bf D}}

\def\size{{\rm size}}
\def\diam{{\rm diam}}
\def\energy{{\rm energy}}
\def\modenergy{\widetilde{\rm energy}}

\def\pv{{\vec P}}
\def\qv{{\vec Q}}
\def\Pv{{\vec{\bf P}}}
\def\Qv{{\vec{\bf Q}}}
\def\dist{{\rm dist}}

\def\111{\gamma}

\def\be#1{\begin{equation}\label{#1}}
\def\bas{\begin{align*}}
\def\eas{\end{align*}}
\def\bi{\begin{itemize}}
\def\ei{\end{itemize}}
\newenvironment{proof}{\noindent {\bf Proof:} }{\endprf\par}
\def \endprf{\hfill  {\vrule height6pt width6pt depth0pt}\medskip}
\def\emph#1{{\it #1}}

\title{$L^p$ estimates for the biest II. The Fourier case}

\author{Camil Muscalu}
\address{Department of Mathematics, UCLA, Los Angeles CA 90095-1555}
\email{camil@@math.ucla.edu}

\author{Terence Tao}
\address{Department of Mathematics, UCLA, Los Angeles CA 90095-1555 }
\email{tao@@math.ucla.edu}

\author{Christoph Thiele}
\address{Department of Mathematics, UCLA, Los Angeles CA 90095-1555}
\email{thiele@@math.ucla.edu}

\include{psfig}
\begin{document}

\begin{abstract}  
We prove $L^p$ estimates (Theorem \ref{main}) for the ``biest'', 
a trilinear multiplier operator with singular symbol.  The methods used 
are based on the treatment of the Walsh analogue of the biest in the 
prequel \cite{mtt:walshbiest} of this paper, but with additional technicalities due to the fact that in the Fourier model one cannot obtain perfect localization in both space and frequency.
\end{abstract}

\maketitle

\section{introduction}

The bilinear Hilbert transform can be written (modulo minor modifications) as
$$
B(f_1,f_2)(x):=\int_{\xi_1<\xi_2}
\widehat{f}_1(\xi_1)\widehat{f}_2(\xi_2) e^{2\pi ix(\xi_1+\xi_2)}\,d\xi_1 d\xi_2,$$
where $f_1, f_2$ are test functions on $\R$ and the Fourier transform is defined by
$$ \hat f(\xi) := \int_\R e^{-2\pi i x \xi} f(x)\ dx.$$
From the work of Lacey and Thiele \cite{laceyt1}, \cite{laceyt2} we have the following $L^p$ estimates on $B$:

\begin{theorem}\label{bht} 
$B$ maps $L^p \times L^q \to L^r$ whenever $1 < p, q \leq \infty$, $1/p + 1/q = 1/r$, and $2/3 < r < \infty$.
\end{theorem}

In this paper we shall study a trilinear variant $T$ of the bilinear Hilbert transform,
defined by 
\begin{equation}\label{oper}
T(f_1,f_2,f_3)(x):=\int_{\xi_1<\xi_2<\xi_3}
\widehat{f}_1(\xi_1)\widehat{f}_2(\xi_2)\widehat{f}_3(\xi_3)
e^{2\pi ix(\xi_1+\xi_2+\xi_3)}\,d\xi_1 d\xi_2 d\xi_3.
\end{equation}
The operator $T$ arises naturally from WKB expansions of eigenfunctions of one-dimensional Schr\"{o}dinger operators, following the work of Christ and Kiselev \cite{ck}; see Appendix I of \cite{mtt:walshbiest} for further discussion.  For these applications it is of interest to obtain $L^p$ estimates on $T$, especially in the case when the functions $f_j$ are in $L^2$.

From the identity
$$
f_1(x) f_2(x) f_3(x) = \int
\widehat{f}_1(\xi_1)\widehat{f}_2(\xi_2)\widehat{f}_3(\xi_3)
e^{2\pi ix(\xi_1+\xi_2+\xi_3)}\,d\xi_1 d\xi_2 d\xi_3
$$
we see that $T$ has the same homogeneity as the pointwise product operator, and hence we expect estimates of H\"older type, i.e. that $T$ maps $L^{p_1} \times L^{p_2} \times L^{p_3}$ to $L^{p'_4}$ when $1/p'_4 = 1/p_1 + 1/p_2 + 1/p_3$.

If the restriction $\xi_1 < \xi_2 < \xi_3$ were replaced by some smoother cutoff then it would be possible to obtain these $L^p$ estimates from the paraproduct theory of Coifman and Meyer (see e.g. \cite{coifman}).  However, as for the bilinear Hilbert transform, the rough nature of the cutoff makes the treatment of this operator much more delicate. 

The operator $T$ is connected to $B$ in several ways.  For instance, we have $T(1, f_2, f_3) = B(P f_2, P f_3)$, where $P$ is the Riesz projection to the half-line $[0,\infty)$.  Note that if we replace the restriction $\xi_1 < \xi_2 < \xi_3$ by just $\xi_1 < \xi_2$ or $\xi_2 < \xi_3$, the operator $T$ factors into $B(f_1,f_2) f_3$ or $f_1 B(f_2,f_3)$.  Thus $T$ can be thought of as a hybrid of the above two operators, although it does not seem to be as easily factorized as the preceding examples.  

The main purpose of this paper is to obtain a large set of $L^p$ estimates for  $T$. As a by-product of our framework we shall be able to give a fairly short (but standard) proof of Theorem \ref{bht}.

Let us consider now the $3$-dimensional affine hyperspace
\[
S:=\{(\alpha_1,\alpha_2,\alpha_3,\alpha_4)\in\R^4\,
|\,\alpha_1 + \alpha_2 + \alpha_3 + \alpha_4=1\}.
\]


Denote by $\D'$ the open interior of the 
convex hull of the $12$ extremal points
$A_1,...,A_{12}$. 
They belong to $S$ and have the following coordinates:

\[
\begin{array}{llll}
A_1:(1,\frac 12,1,-\frac 32)  &  A_2:(\frac 12,1,1,-\frac 32)  &  
A_3:(\frac 12,1,-\frac 32,1) & A_4:(1,\frac 12,-\frac 32,1)  \\  
\ &\ &\ & \ \\
A_5:(1,-\frac 12,0,\frac 12)  &  A_6:(1,-\frac 12,\frac 12,0) & 
A_7:(\frac 12,-\frac 12,0,1)  &  A_8:(\frac 12,-\frac 12,1,0)  \\  
\ &\ &\ & \ \\
A_9:(-\frac 12,1,0,\frac 12) &
A_{10}:(-\frac 12,1,\frac 12,0) &  A_{11}:(-\frac 12,\frac 12,1,0) &  A_{12}:(-\frac 12,\frac 12,0,1).
\end{array}
\]

See \cite{mtt:walshbiest} for a
figure of this polytope.


Also, denote by $\D''$ the open interior of the 
convex hull of the $12$ extremal points
$B_1,...,B_{12}$ in $S$ (they are not represented in the picture)
where the coordinates of $B_j$ are obtained from the coordinates of
$A_j$ after permuting the indices $1$ and $3$ for $j=1,...,12$ (for
instance $B_2$ has the coordinates $B_2(1,1,1/2,-3/2)$).
Then set $\D:= \D'\cap \D''$.

Our main result is:

\begin{theorem}\label{main}
Let $1<p_1, p_2, p_3\leq\infty$ and $0<p'_4<\infty$ such that
\[\frac{1}{p_1}+\frac{1}{p_2}+\frac{1}{p_3}=\frac{1}{p'_4}.\]
Then $T$ maps
\begin{equation}\label{star}
T: L^{p_1}\times L^{p_2}\times L^{p_3}\rightarrow L^{p'_4}
\end{equation}
as long as $(1/p_1,1/p_2,1/p_3,1/p_4)\in \D$.
\end{theorem}
In particular,
$T$ maps $L^{p_1}\times L^{p_2}\times L^{p_3}\rightarrow L^{p'_4}$,
whenever $1<p_1, p_2, p_3\leq\infty$ and $1\leq p'_4 <\infty$.

For the particular application to Schr\"odinger eigenfunctions, the functions $f_i$ will be in $L^2$.  We thus record the corollary 

\begin{corollary}\label{23}
$T$ maps $L^2\times L^2\times L^2\rightarrow L^{2/3}$.
\end{corollary}

In \cite{mtt:walshbiest} Theorem \ref{main} was proved for a Walsh-Fourier analogue $T_{walsh,\P,\Q}$ of $T$.  From the point of view of time-frequency analysis the two operators are closely related, however the Walsh model is easier to analyze technically because it is possible in that model to localize perfectly in both space and frequency simultaneously.  In the Fourier case one has to deal with several ``Schwartz tails'' which introduce additional difficulties.  For instance, in the Walsh model an inner product $\langle \phi_P, \phi_Q \rangle$ of wave packets vanishes unless the spatial intervals $I_P$ and $I_Q$ are nested; however in the Fourier model one needs to consider the case when $I_P$ and $I_Q$ are separated (although the estimates improve rapidly with the relative separation of $I_P$ and $I_Q$).

If $m$ is a function on $\R^3$, then we define the multiplier operator
\[
T_m(f_1,f_2,f_3)(x):=\int m(\xi_1,\xi_2,\xi_3)
\widehat{f}_1(\xi_1)\widehat{f}_2(\xi_2)\widehat{f}_3(\xi_3)
e^{2\pi ix(\xi_1+\xi_2+\xi_3)}\,d\xi_1 d\xi_2 d\xi_3.
\]
In \cite{cct} such operators were studied, in particular
estimates as in Theorem \ref{main} were proved if $m$ satisfies
\begin{equation}\label{symb-est}
|\partial_\xi^\alpha m(\xi)| 
\lesssim |\dist(\xi, \{\lambda,\lambda,\lambda\})|^{-|\alpha|}
\end{equation}
for all multi-indices $\alpha$.

As in \cite{mtt:walshbiest} it shall be convenient to split $T$ into two pieces, plus an error term,
which is a multiplier operator of the type \eqref{symb-est}.  Specifically, we let $m'(\xi)$ 
be a function which coincides with the characteristic function
$\chi_{\xi_2<\xi_3}$ on a cone 
\[\{|\xi_3-\xi_2|\ll |\xi_1-\frac{\xi_2+\xi_3}2|, \xi_1<\xi_2\}\]
and satisfies \eqref{symb-est} outside a smaller cone of the same type
(different implicit constant).
Likewise, we let $m''(\xi)$ be a function which coincides with the
characteristic function $\chi_{\xi_1<\xi_2}$ on a cone 
\[\{|\xi_2-\xi_1|\ll |\xi_3-\frac{\xi_1+\xi_2}2|, \xi_2<\xi_3\}\]
and satisfies (\ref{symb-est}) outside a smaller cone of the same type.
Then
\[\chi_{\xi_1<\xi_2<\xi_3}-m'-m''\]
satisfies \eqref{symb-est} everywhere and can be estimated\footnote{One can also bound $T'''$ by the arguments in this paper, indeed the treatment would be simpler as the tri-tiles $\Pv$ and $\Qv$ will be constrained to be roughly the same size.} by
the results in {\cite{cct}}. The multipliers $m'$ and $m''$ will be specified later.

Thus to prove our main theorem it suffices to show

\begin{theorem}\label{main-p}
Let $1<p_1, p_2, p_3\leq\infty$ and $0<p'_4<\infty$ such that
\[\frac{1}{p_1}+\frac{1}{p_2}+\frac{1}{p_3}=\frac{1}{p'_4}.\]
Then $T_{m'}$ maps
\begin{equation}\label{star-p}
T_{m'}: L^{p_1}\times L^{p_2}\times L^{p_3}\rightarrow L^{p'_4}
\end{equation}
as long as $(1/p_1,1/p_2,1/p_3,1/p_4)\in \D'$.  Similarly for $T_{m''}$ and $\D''$.
\end{theorem}

We shall only prove the claim for $T_{m'}$, as the claim for $T_{m''}$ 
follows by a  permutation of the 1 and 3 indices.

We remark that Theorem \eqref{main} is about a trilinear multiplier operator whose
symbol is the characteristic function of the set $\xi_1<\xi_2<\xi_3$.
Modifications as in \cite{gilbertnahmod} give theorems for more
general multiplier symbols. However, a product structure condition
as discussed in \cite{grafakoskalton} seems to occur. While these are questions
of interest, we shall not discuss them further in this article.

There are various recurring themes in the subject of multilinear 
singular integrals as in 
\cite{gilbertnahmod}, \cite{grafakosli},
\cite{lacey}, \cite{laceyt1}, \cite{laceyt2}, \cite{cct}
\cite{thiele'}, \cite{thiele}
and so forth, which the current paper 
again builds up on. While the current article is mostly self contained,
we will mark as ``standard'' any arguments that are well 
understood by now in this framework.  Our notation and setup follows
closely \cite{cct} and \cite{mtt:walshbiest}.

The authors would like to thank Mike Christ for pointing out to them 
the occurence of multilinear singular integrals of the type discussed
in this article in the study of eigenfunction expansions of 
Schr\"odinger operators.

The first author was partially supported by a Sloan Dissertation Fellowship.
The second author is a Clay Prize Fellow and is supported by grants from
the Sloan and Packard Foundations. The third author was partially supported
by a Sloan Fellowship and by NSF grants DMS 9985572 and DMS 9970469.

\section{Notation}\label{notation-sec}

In this section we set out some general notation used throughout the paper.  

We use $A\lesssim B$ to denote the statement that
$A\leq CB$ for some large constant $C$, and $A\ll B$ to denote the
statement that $A\leq C^{-1}B$ for some large constant $C$. Our constants $C$ shall always be independent of $\P$ and $\Q$.  Given any interval (or cube) $I$, we let $|I|$ denote the measure of $I$, and $cI$ denotes the interval (or cube) with the same center as $I$ but $c$ times the side-length.

Given a spatial interval $I$, we shall define the approximate cutoff function $\tilde \chi_I$ by
$$ \tilde \chi_I(x) := (1 + (\frac{|x-x_I|}{|I|})^2)^{-1/2},$$
where $x_I$ is the center of $I$.

A collection $\{\omega\}$ of intervals is said to be \emph{lacunary around the frequency $\xi$} if we have  $\dist(\xi, \omega) \sim |\omega|$ for all $\omega$ in the collection.

We define a \emph{modulated Calder\'on-Zygmund operator} to be any operator $T$ which is bounded on $L^2$ and has the form
$$ Tf(x) = \int K(x,y) f(y)\ dy$$
where $x,y \in \R$, and the (possibly vector-valued) kernel $K$ obeys the estimates 
$$ |K(x,y)| \lesssim 1/|x-y|$$
and
$$ |\nabla_{x,y} (e^{2\pi i (x\xi + y\eta)} K(x,y))| \lesssim 1/|x-y|^2$$
for all $x \neq y$ and for some $\xi,\eta \in \R$.  Note that a modulated Calder\'on-Zygmund operator is the composition of an ordinary Calder\'on-Zygmund operator with modulation operators such as $f \mapsto e^{2\pi i \xi \cdot} f$.  By standard Calder\'on-Zygmund theory (see e.g. \cite{stein}) we thus see that $T$ is bounded on $L^p$ for all $1 < p < \infty$, and is also weak-type $(1,1)$.

\section{interpolation}\label{interp-sec}

In this section we review the interpolation theory from \cite{cct} which allows us to reduce multi-linear $L^p$ estimates such as those in Theorem \ref{main} to certain ``restricted weak type'' estimates.  This will be the first reduction in our proof of Theorem \ref{main}.

To prove the $L^p$ estimates on $T'$ it is convenient to use duality and introduce the quadrilinear form $\Lambda'$ associated to $T'$ via the formula
$$
\Lambda'(f_1, f_2,f_3, f_4) := \int_{\R} T'(f_1, f_2, f_3)(x)f_4(x) dx.
$$
The statement that $T'$ is bounded from $L^{p_1} \times L^{p_2} \times L^{p_3}$ to $L^{p'_4}$ is then equivalent to $\Lambda'$ being bounded on $L^{p_1} \times L^{p_2} \times L^{p_3} \times L^{p_4}$ if $1 < {p_4}' < \infty$.  
For ${p_4}' <1$ this simple duality relationship breaks down, however the interpolation arguments in \cite{cct} will allow us to reduce (\ref{star}) to certain ``restricted
type'' estimates on $\Lambda'$. As in \cite{cct} we find more
convenient to work with the quantities $\alpha_i=1/p_i$, $i=1,2,3,4$,
where $p_i$ stands for the exponent of $L^{p_i}$.

\begin{definition}
A tuple $\alpha=(\alpha_1,\alpha_2,\alpha_3,\alpha_4)$ is called admissible, if

\[-\infty <\alpha_i <1\]
for all $1\leq i\leq 4$,

\[\sum_{i=1}^4\alpha_i=1\]
and there is at most one index $j$ such that $\alpha_j<0$. We call
an index $i$ good if $\alpha_i\geq 0$, and we call it bad if
$\alpha_i<0$. A good tuple is an admissible tuple without bad index, a
bad tuple is an admissible tuple with a bad index.
\end{definition}

\begin{definition}
Let $E$, $E'$ be sets of finite measure. We say that $E'$ is a major 
subset of $E$ if $E'\subseteq E$ and $|E'|\geq\frac{1}{2}|E|$.
\end{definition} 

\begin{definition}
If $E$ is a set of finite measure, we denote by $X(E)$ the space of
all functions $f$ supported on $E$ and such that $\|f\|_{\infty}\leq
1$.
\end{definition}

\begin{definition}
If $\alpha=(\alpha_1,\alpha_2,\alpha_3,\alpha_4)$ is an admissible
tuple, we say
that a $4$-linear form $\Lambda$ is of restricted type $\alpha$
if for every sequence $E_1, E_2, E_3, E_4$ of subsets of $\R$ with
finite measure, there exists a major subset $E'_j$ of $E_j$ for each
bad index $j$ (one or none) such that

\[|\Lambda(f_1, f_2, f_3, f_4)|\lesssim |E|^{\alpha}\]
for all functions $f_i\in X(E'_i)$, $i=1,2,3,4$, where we adopt the
convention $E'_i =E_i$ for good indices $i$, and $|E|^{\alpha}$
is a shorthand for

\[|E|^{\alpha}=|E_1|^{\alpha_1}|E_2|^{\alpha_2}
|E_3|^{\alpha_3}|E_4|^{\alpha_4}.\]
\end{definition}

The following ``restricted type'' result will be proved directly.

\begin{theorem}\label{teorema1}
For every vertex $A_i$, $i = 1, \ldots, 12$ there exist
admissible tuples $\alpha$ arbitrarily close to $A_i$ such that the form $\Lambda'$ is of restricted
type $\alpha$.
\end{theorem}

By interpolation of restricted weak type estimates (cf. \cite{cct}) we thus obtain

\begin{corollary}\label{restricted}
Let $\alpha$ be an admissible tuple such that $\alpha\in \D'$. 
Then $\Lambda'$ is of restricted type
$\alpha$.
\end{corollary}

It only remains to convert these restricted type estimates into strong type estimates \eqref{star-p}.  To do this one just has to apply (exactly as in \cite{cct}) the multilinear Marcinkiewicz interpolation theorem \cite{janson} in the case of good tuples and the interpolation lemma 3.11 in \cite{cct} in the case of bad tuples.

This ends the proof of theorem \ref{main-p}. Hence, it remains to prove
Theorem \ref{teorema1}.

\section{Discretization}\label{discrete-sec}

We now prove Theorem \ref{teorema1}.  The first reduction is to pass from the ``continuous'' form $\Lambda'$ to a ``discretized'' variant involving sums of inner products with wave packets.  This step is standard 
and appears essential in order for the phase plane combinatorics to work correctly.  

\begin{definition}\label{mesh-def}
Let $n \geq 1$ and $\sigma \in \{ 0, \frac{1}{3}, \frac{2}{3} \}^n$.  We define the \emph{shifted $n$-dyadic mesh} $D = D^n_{\sigma}$ to be the collection of cubes of the form
$$ D^n_{\alpha} := \{ 2^j (k + (0,1)^n + (-1)^j \sigma)| j \in \Z, \quad k \in \Z^n \}.$$
We define a \emph{shifted dyadic cube} to be any member of a shifted $n$-dyadic mesh.
\end{definition}

Observe that for every cube $Q$, there exists a shifted dyadic cube $Q'$ such that $Q \subseteq \frac{7}{10} Q'$ and $|Q'| \sim |Q|$; this is best seen by first verifying the $n=1$ case. 

\begin{definition}\label{sparse-grid}
A subset $D'$ of a shifted $n$-dyadic grid $D$ is called sparse, if for any two 
cubes $Q,Q'$ in $D$ with $Q\neq Q'$ we have $|Q|<|Q'|$ implies $|10^{9}Q|<|Q'|$ and 
$|Q|=|Q'|$ implies $10^9Q\cap 10^9 Q'=\emptyset$.
\end{definition}

Observe that any subset of a shifted $n$-dyadic grid (with $n\le 3$ say), can be 
split into $O(1)$ sparse subsets.

\begin{definition}  Let $\sigma = (\sigma_1,\sigma_2,\sigma_3) \in \{ 0, \frac{1}{3}, \frac{2}{3} \}^3$, and let $1 \leq i \leq 3$.  An $i$-tile with shift $\sigma_i$ is a rectangle $P = I_P \times \omega_P$ with area 1 and with $I_P \in D^1_0$, $\omega_P \in D^1_{\sigma_i}$.  A tri-tile with shift $\sigma$ is a $3$-tuple $\pv = (P_1, \ldots, P_3)$ such that each $P_i$ is an $i$-tile with shift $\sigma_i$, and the $I_{P_i} = I_\pv$ are independent 
of $i$.  The frequency cube $Q_\pv$ of a tri-tile is defined to be
$\prod_{i=1}^3 \omega_{P_i}$.
\end{definition}

We shall sometimes refer to $i$-tiles with shift $\sigma$ just as $i$-tiles, or even as tiles, if the parameters $\sigma$, $i$ are unimportant.

\begin{definition} A set $\Pv$ of tri-tiles is called sparse, if all tri-tiles in $\Pv$
have the same shift and the set $\{Q_{\pv}:\pv\in \Pv\}$ is sparse.
\end{definition}

Again, any set of tri-tiles can be split into $O(1)$ sparse subsets.

\begin{definition}  Let $P$ and $P'$ be tiles.  We write $P' < P$ if $I_{P'} \subsetneq I_P$ and $3\omega_P \subseteq 3\omega_{P'}$, and $P' \leq P$ if $P' < P$ or $P' = P$.
We write $P' \lesssim P$ if $I_{P'} \subseteq I_P$ and $10^7\omega_P \subseteq 10^7 \omega_{P'}$.  We write $P' \lesssim' P$ if $P' \lesssim P$ and
$P' \not \leq P$.
\end{definition}

The ordering $<$ is in the spirit of that in
Fefferman \cite{fefferman} or Lacey and Thiele \cite{laceyt1}, \cite{laceyt2}, \cite{thiele}, but slightly different as $P'$ and $P$ do not quite have to intersect.  This is more convenient for technical purposes.

\begin{definition}\label{rank-def}  A collection $\Pv$ of tri-tiles is said to have \emph{rank 1} if one has the following properties for all $\pv, \pv' \in \Pv$:
\begin{itemize}
\item If $\pv \neq \pv'$, then $P_j \neq P'_j$ for all $j=1,2,3$.
\item If $P'_j \leq P_j$ for some $j=1,2,3$, then $P'_i \lesssim P_i$ for all $1 \leq i \leq 3$.  
\item If we further assume that $|I_{\pv'}| < 10^9 |I_{\pv}|$, then we have $P'_i \lesssim' P_i$ for all $i \neq j$.
\end{itemize}
\end{definition}

\begin{definition}\label{packet-def}  Let $P$ be a tile.  A \emph{wave packet on $P$} is a function $\phi_P$ which has Fourier support in $\frac{9}{10} \omega_P$ and obeys the estimates 
\be{decay}
|\phi_P(x)| \lesssim |I_P|^{-1/2} \tilde \chi_I(x)^M
\end{equation}
for all $M > 0$, with the implicit constant depending on $M$.
\end{definition}

Heuristically, $\phi_P$ is $L^2$-normalized and is supported in $P$.

The discretized form of Theorem \ref{teorema1} is as follows.

\begin{theorem}\label{discrete}
Let $\sigma, \sigma' \in \{ 0, \frac{1}{3}, \frac{2}{3} \}^3$ be shifts, and let $\Pv$, $\Qv$ be finite collections of multi-tiles with shifts $\sigma$, $\sigma'$ respectively such that $\Pv$ and $\Qv$ both have rank 1.  For each $i=1,2,3$ and $\pv \in \Pv$, let $\phi_{P_i} = \phi_{P_i,i}$ be a wave packet on $P_i$.  Similarly for each $i=1,2,3$ and $\qv \in \Qv$ let $\tilde\phi_{Q_i} = \tilde\phi_{Q_i,i}$ be a wave packet on $Q_i$.  Define the form $\Lambda_{\Pv,\Qv}$ by
$$\Lambda_{\Pv,\Qv}(f_1,f_2,f_3,f_4) :=
\sum_{\pv \in \Pv} \frac{1}{|I_\pv|^{1/2}} 
\langle f_1, \phi_{P_1} \rangle
\langle B_{P_2}(f_2,f_3), \phi_{P_2} \rangle
\langle f_4, \phi_{P_3} \rangle
$$
where
$$
B_{P_2}(f_2,f_3) := \sum_{\qv \in \Qv: \omega_{Q_3} \subseteq \omega_{P_2}}
\frac{1}{|I_\qv|^{1/2}}
\langle f_2, \tilde \phi_{Q_1} \rangle
\langle f_3, \tilde \phi_{Q_2} \rangle
\tilde \phi_{Q_3}.$$
Then $\Lambda$ is of restricted type
$\alpha$ for all admissible tuples $\alpha \in \D'$, uniformly in the parameters $\sigma$, $\sigma'$, $\Pv$, $\Qv$, $\phi_{P_i}$, $\phi'_{Q_i}$. Furthermore, in the case that $\alpha$ has a bad index $j$, the restricted type is uniform in the sense that the  major subset $E'_j$ can be chosen independently of the parameters just mentioned.
\end{theorem}

This theorem is the Fourier analogue of Theorem 2.5 of \cite{mtt:walshbiest}.  The main new difficulty is that unlike the Walsh wave packets $\phi_{walsh,P}$, the Fourier wave packets $\phi_P$ are not perfectly localized in physical space to the interval $I_P$.

In the rest of this section we show how Theorem \ref{teorema1} can be deduced from Theorem \ref{discrete}.  This section can be read independently of the later sections, which are concerned with the proof of Theorem \ref{discrete}.

Assume Theorem \ref{discrete} holds.  This Theorem is phrased in terms of wave packets which are perfectly localized in frequency and imperfectly localized in physical space.  We now bootstrap this theorem to a similar statement in which no physical space localization is made.

\begin{definition}\label{cube-rank-def}  Let $\sigma \in \{ 0, \frac{1}{3}, \frac{2}{3} \}^3$ be a shift.  A collection $\Q \subset D^3_\sigma$ of cubes is said to have \emph{rank 1} if one has the following properties for all $Q, Q' \in \Q$:
\begin{itemize}
\item If $Q \neq Q'$, then $Q \cap Q' = \emptyset$.  (In other words, the $Q$ are disjoint).
\item If $Q \neq Q'$, then $Q_i \neq Q'_i$ for all $i=1,2,3$.
\item If $3Q'_j \subset 3Q_j$ for some $j=1,2,3$, then $10^7Q'_i \subset 10^7Q_i$ for all $1 \leq i \leq 3$.  
\item If we further assume that $|Q'| < |10^9Q|$, then we have $3Q'_i \cap 3Q_i = \emptyset$ for all $i \neq j$.
\end{itemize}
\end{definition}

If $Q$ is a cube in $\R^3$, denote by $-Q$ the reflected 
cube about the origin; similarly for intervals in $\R$.

\begin{corollary}\label{discrete-cor}
Let $\sigma, \sigma' \in \{ 0, \frac{1}{3}, \frac{2}{3} \}^3$ be shifts, let $\Q$, $\Q'$ be finite rank 1 collections of cubes in $D^3_\sigma$, $D^3_{\sigma'}$.  For each $i=1,2,3$ and $Q \in \Q$, let $\eta_{Q_i} = \eta_{Q_i,i}$ be a bump function adapted to $\frac{9}{10}Q_i$.  Similarly for each $i=1,2,3$ and $Q' \in \Q'$ let $\eta'_{Q'_i} = \eta'_{Q'_i,i}$ be a bump function adapted to $\frac{9}{10} Q'_i$.  Then the form
$$\Lambda(f_1,f_2,f_3,f_4) :=
\int \delta(\xi_1+\xi_2+\xi_3+\xi_4) m(\xi_1,\xi_2,\xi_3,\xi_4)
\widehat{f}_1(\xi_1)\widehat{f}_2(\xi_2)\widehat{f}_3(\xi_3)
\widehat{f}_4(\xi_4)\, d\xi_1,\dots, d\xi_4$$
is of restricted type
$\alpha$ for all admissible tuples $\alpha \in \D'$, uniformly in the parameters $\sigma$, $\sigma'$, $\Q$, $\Q'$, $\eta_{Q_i}$, $\eta'_{Q'_i}$, where
$$ m(\xi_1,\xi_2,\xi_3,\xi_4) :=
\sum_{Q \in \Q} \sum_{Q' \in \Q': -Q'_3 \subseteq Q_2}
\eta_{Q_1}(\xi_1) \eta_{Q_2}(\xi_2+\xi_3) \eta_{Q_3}(\xi_4)
\eta'_{Q'_1}(\xi_2) \eta'_{Q'_2}(\xi_3) \eta'_{Q'_3}(\xi_1 + \xi_4).$$
Furthermore, in the case that $\alpha$ has a bad index $j$, the restricted type is uniform in the sense that the  major subset $E'_j$ can be chosen independently of the parameters just mentioned.
\end{corollary}

\begin{proof}
To motivate matters, let us first consider the simpler object
$$
\left<f_2,a_Q\right>
:=\int \delta(\xi_1+\xi_2+\xi_3) 
\eta_{Q_1}(\xi_1) \eta_{Q_2}(\xi_2) \eta_{Q_3}(\xi_3)
\widehat{f}_1(\xi_1)\widehat{f}_2(\xi_2)\widehat{f}_3(\xi_3)
d\xi_1d\xi_2d\xi_3$$
for some cube $Q \in \Q$.  By Plancherel this is equal to
\bas
&
\int (f_1 * \check \eta_{Q_1})(x)
(f_2 * \check \eta_{Q_2})(x)
(f_3 * \check \eta_{Q_3})(x)
\ dx\\
&= l(Q)^{3/2}
\int \langle f_1, \phi_{Q_1,x} \rangle
\langle f_2, \phi_{Q_2,x} \rangle
\langle f_3, \phi_{Q_3,x} \rangle\ dx
\end{align*}
where
$$ \phi_{Q_j,x}(y) := l(Q)^{-1/2} \overline{\check \eta_{Q_j}(x-y)}$$
and $l(Q)$ is the side-length of $Q$.  We can rewrite this as
$$
\int_0^1  \sum_{\pv: Q_\pv = Q} |I_\pv|^{-1/2}
\langle f_1, \phi_{P_1,t,1} \rangle
\langle f_2, \phi_{P_2,t,2} \rangle
\langle f_3, \phi_{P_3,t,3} \rangle
\ dt
$$
where $\pv$ ranges over all tri-tiles with frequency cube $Q$ and spatial interval $I_\pv$ in $D^1_0$, $\phi_{P_j,t,j}$ is the function
$$ \phi_{P_j,t,j} := \phi_{Q_j,x_{\pv} + |I_\pv| t}$$
and $x_{\pv}$ is the center of $I_\pv$.  Note that $\phi_{P_j,t,j}$ is a wave packet on $P_j$ uniformly in $t$.

Similarly, consider
$$\left<f_3,b_{Q'}\right>
:=\int \delta(\xi_1+\xi_2+\xi_3) 
\eta_{Q_1'}(\xi_1) \eta_{Q_2'}(\xi_2) \eta_{Q_3'}(\xi_3)
\widehat{f}_1(\xi_1)\widehat{f}_2(\xi_2)\widehat{f}_3(\xi_3)
d\xi_1d\xi_2d\xi_3$$

The multiplier $m$ in the statement of the Corollary can now be rewritten as
$$
\Lambda(f_1,f_2,f_3,f_4)=\sum_{Q\in \Q,Q'\in \Q':-Q_3'\subseteq Q_2}
\int \overline{\widehat{a_Q}(\tau)\widehat{b_{Q'}}(-\tau)} \, d\tau
$$
$$
=\int_0^1 \int_0^1 \sum_{\pv: Q_\pv \in \Q}\frac{1}{|I_\pv|^{1/2}} 
\langle f_1, \phi_{P_1,t,1} \rangle
\langle B_{P_2,t'}(f_2,f_3), \phi_{P_2,t,2} \rangle
\langle f_4, \phi_{P_3,t,3} \rangle\ dt dt'$$
where
$$B_{P_2,t'}(f_2,f_3) := \sum_{\qv: Q_\qv \in \Q', - \omega_{Q_3} \subseteq \omega_{P_2}}
\frac{1}{|I_\qv|^{1/2}}
\langle f_2, \phi_{Q_1,t',1} \rangle
\langle f_3, \phi_{Q_2,t',2} \rangle
\overline{\phi_{Q_3,t',3}}.$$
Note that the collection of tri-tiles $\pv$ has rank 1, and similarly for the collection of tri-tiles $\qv$.  Observe that we can get rid of the complex conjugation sign in the definition of $B_{P_2,t'}$ by redefining 
$Q_3$ to be $-Q_3$ and 
redefining $\phi_{Q_3,t',3}$ accordingly; 
this 
also replaces the condition
$- \omega_{Q_3} \subseteq \omega_{P_2}$ by the condition $\omega_{Q_3} \subseteq \omega_{P_2}$, but it does not change the rank one property of the collection $\Qv$.
The claim then follows by integrating the conclusion of Theorem \ref{discrete} over $t$, $t'$, using the uniformity assumptions of that Theorem.  (The finiteness condition on $\Pv$ and $\Qv$ can be removed by the usual limiting arguments.)
\end{proof}


We can now prove Theorem \ref{teorema1}.  By a standard partition of unity we can write
$$
\chi_{2\xi_1 < \xi_2}(\xi_1,\xi_2,\xi_3) = 
\sum_{\sigma \in \{0, \frac 1 3, \frac 2 3\}^3} \sum_{Q \in \Q_\sigma} \phi_{Q,\sigma}(\xi_1, \xi_2, \xi_3)$$
whenever $\xi_1 + \xi_2 + \xi_3 = 0$, where $\Q_\sigma \subset D^3_\sigma$ is a collection of cubes which intersect the plane $\{\xi_1 + \xi_2 + \xi_3 = 0\}$ and which satisfy the Whitney property
$$ 10^3 \diam(Q) \leq \dist(Q, \{ 2\xi_1 = \xi_2, \xi_1 + \xi_2 + \xi_3 = 0 \}) \leq 10^5 \diam(Q) $$
for all $Q \in \Q_\sigma$, and for each cube $Q \in \Q_\sigma$, $\phi_{Q,\sigma}$ is a bump function adapted to $\frac{8}{10} Q$.  Note that by refining $\Q_\sigma$ by a finite factor if necessary one can make $\Q_\sigma$ have rank 1.

By splitting $\phi_{Q,\sigma}$ as a Fourier series in the $\xi_i$ we can then write
$$
\chi_{2\xi_1 < \xi_2}(\xi_1,\xi_2,\xi_3) = \sum_{k \in \Z^3} c_k \sum_{\sigma \in \{0, \frac 1 3, \frac 2 3\}^3} \sum_{Q \in \Q_\sigma} \eta_{Q_1,\sigma,k,1}(\xi_1)
\eta_{Q_2,\sigma,k,2}(\xi_2)
\eta_{Q_3,\sigma,k,3}(\xi_3),$$
where $c_k$ is a rapidly decreasing sequence 
and $\eta_{Q_j, \alpha, k, j}$ is a bump function adapted to $\frac{9}{10} Q_j$ uniformly in $k$  
(the $c_k$ are fractional powers of the Fourier coefficients, 
whereas the $\eta_{Q_j,\sigma,k,j}$ are products
of an initial bump function, times a complex exponential, times another
fractional power of the Fourier coefficient which generates the 
uniformity of the bump function in $k$).

Similarly, we can write
$$
\chi_{\xi_1 < \xi_2}(\xi_1,\xi_2,\xi_3) = \sum_{k \in \Z^3} c'_k \sum_{\sigma \in \{0, \frac 1 3, \frac 2 3\}^3} \sum_{Q \in \Q_\sigma} \eta'_{Q_1,\sigma,k,1}(\xi_1)
\eta'_{Q_2,\sigma,k,2}(\xi_2)
\eta'_{Q_3,\sigma,k,3}(\xi_3),$$

As a consequence, the expression 
\bas
\sum_{k,k' \in \Z^3} c_k c'_{k'} \sum_{\sigma,\sigma' \in \{0, \frac 1 3, \frac 2 3\}^3} \sum_{Q \in \Q_\sigma} \sum_{Q' \in \Q_{\sigma'}: -Q'_3\subseteq Q_2} &\eta_{Q_1,\alpha,k,1}(\xi_1)
\eta_{Q_2,\alpha,k,2}({\xi_2+\xi_3})
\eta_{Q_3,\alpha,k,3}(\xi_4)\\
&\eta_{Q_1',\alpha,k',1}(\xi_2)
\eta_{Q_2',\alpha,k',2}(\xi_3)
\eta_{Q_3',\alpha,k',3}(\xi_1+\xi_4)
\end{align*}
on the hyperplane $\xi_1 + \xi_2 + \xi_3 + \xi_4=0$ is equal to $\chi_{2\xi_1 < \xi_2 + \xi_3} \chi_{\xi_2 < \xi_3}$ when $|\xi_3 - \xi_2| \ll |\xi_1 - \frac{\xi_2 + \xi_3}{2}|$, (under the latter constraint the condition $-Q_3'\subseteq Q_2$ is automatic for nonzero summands), 
and satisfies \eqref{symb-est} outside any cone of the type
$|\xi_3 - \xi_2| \ll |\xi_1 - \frac{\xi_2 + \xi_3}{2}|$ (then the scales of the
cubes $Q$ and $Q'$ are essentially coupled and also coupled to
the distance to the line $\xi_1=\xi_2=\xi_3$.)

Thus the above expression is a multiplier $m'(\xi_1, \xi_2, \xi_3, \xi_4)$ of the type requested for the definition of $T'$.  Theorem \ref{teorema1} then follows from Corollary \ref{discrete-cor} and summing in the parameters $k,k',\sigma,\sigma'$, using the uniformity conclusions in Corollary \ref{discrete-cor}.  (The finiteness assumption on $\Q$ and $\Q'$ in Corollary \ref{discrete-cor} can be removed by the usual limiting arguments).

\section{trees}

The standard approach to prove the desired estimates for the forms $\Lambda_{\Pv,\Qv}$
is to organize our collections of tri-tiles $\Pv$, $\Qv$ into trees as in \cite{fefferman}.
We may assume and shall do so for the rest of this article that $\Pv$ and $\Q_v$ are sparse.

\begin{definition} For any $1 \leq j \leq 3$ and a tri-tile $\pv_T \in \Pv$, define a $j$-tree with top $\pv_T$ to be a collection of tri-tiles $T \subseteq \Pv$ such that
$$ P_j \leq P_{T,j} \hbox{ for all } \pv \in T,$$
where $P_{T,j}$ is the $j$ component of $\pv_T$.  We write $I_T$ and $\omega_{T,j}$ for $I_{\pv_T}$ and $\omega_{P_{T,j}}$ respectively.  We say that $T$ is a tree if it is a $j$-tree for some $1 \leq j \leq 3$.
\end{definition}

Note that $T$ does not necessarily have to contain its top $\pv_T$.

\begin{definition} Let $1 \leq i \leq 3$.  Two trees $T$, $T'$ are said to be \emph{strongly $i$-disjoint} if 
\bi
\item $P_i \neq P'_i$ for all $\pv \in T$, $\pv' \in T'$.
\item Whenever $\pv \in T$, $\pv' \in T'$ are such that
$2\omega_{P_i} \cap 2\omega_{P'_i}\neq \emptyset$, then one has
$I_{\pv'} \cap I_T = \emptyset$, and similarly with $T$ and $T'$ reversed.
\end{itemize}
\end{definition}

Note that if $T$ and $T'$ are strongly $i$-disjoint, 
then $I_P\times 2\omega_{P_i} \cap I_{P'}\times 2\omega_{{P'}_i} = \emptyset$ 
for all $\pv \in T$, $\pv' \in T'$.

Given that $\P_v$ is sparse, it is easy to see that
if $T$ is an $i$- tree, then for all $\pv,\pv'\in T$ and $j\neq i$ we have
$$ \omega_{P_j} = \omega_{P'_j}$$ or 
$$ 2\omega_{P_j} \cap 2\omega_{P'_j}=\emptyset$$.

\section{Tile norms}

In the sequel we shall be frequently estimating expressions of the form
\be{trilinear}
| \sum_{\pv \in \Pv} \frac{1}{|I_\pv|^{1/2}} a^{(1)}_{P_1} a^{(2)}_{P_2} a^{(3)}_{P_3}|
\end{equation}
where $\Pv$ is a collection of tri-tiles and $a^{(j)}_{P_j}$ are complex numbers for $\pv \in \Pv$ and $j = 1,2,3$.  In some cases (e.g. if one only wished to treat the Bilinear Hilbert transform) we just have
\be{aj-def}
a^{(j)}_{P_j} = \langle f_j, \phi_{P_j} \rangle
\end{equation}
 but we will have more sophisticated sequences $a^{(j)}_{P_j}$ when dealing with $\Lambda_{\Pv,\Qv}$.  

In \cite{mtt:walshbiest} the following (standard) norms on sequences of tiles were introduced: 

\begin{definition}\label{size-def}
Let $\Pv$ be a finite collection of tri-tiles, $j=1,2,3$, and let $(a_{P_j})_{\pv \in \Pv}$ be a sequence of complex numbers.  We define the \emph{size} of this sequence by
$$ \size_j( (a_{P_j})_{\pv \in \Pv} ) := \sup_{T \subset \Pv}
(\frac{1}{|I_T|} \sum_{\pv \in T} |a_{P_j}|^2)^{1/2}$$
where $T$ ranges over all trees in $\Pv$ which are $i$-trees for some $i \neq j$.
We also define the \emph{energy} of the sequence by
$$ \energy_j((a_{P_j})_{\pv \in \Pv} ) := \sup_{\D \subset \Pv}
(\sum_{\pv \in \D} |a_{P_j}|^2)^{1/2}$$
where $\D$ ranges over all subsets of $\Pv$ such that the tiles $\{ P_j: \pv \in \D \}$ are pairwise disjoint.
\end{definition}

The size measures the extent to which the sequence $a_{P_j}$ can concentrate on a single tree and should be thought of as a phase-space variant of the BMO norm.  The energy is a phase-space variant of the $L^2$ norm.  As the  notation suggests, the number $a_{P_j}$ should be thought of as being associated with the tile $P_j$ rather than the full tri-tile $\pv$.

In the Walsh model the energy is a tractable quantity; for instance, if $a_{P_j}$ is given by \eqref{aj-def} then one can control the energy by $\|f_j\|_2$ thanks to the perfect orthogonality of the Walsh wave packets.  However, in the Fourier case the orthogonality is too poor to give a usable bound on the energy, and so we must instead use a more technical substitute.

\begin{definition}\label{energy-def}
Let the notation be as in Definition \ref{size-def}.  We define the \emph{modified energy} of the sequence $(a_{P_j})_{\pv \in \Pv}$ by
\be{energy-mod-def}
\modenergy_j((a_{P_j})_{\pv \in \Pv} ) := \sup_{n \in \Z}
\sup_{\T}
2^{n} (\sum_{T \in \T} |I_T|)^{1/2}
\end{equation}
where $\T$ ranges over all collections of strongly $j$-disjoint trees in $\Pv$ such that
$$ (\sum_{\pv \in T} |a_{P_j}|^2)^{1/2} \ge 2^n |I_T|^{1/2}$$
for all $T \in \T$, and
$$ (\sum_{\pv \in {T'}} |a_{P_j}|^2)^{1/2} \le 2^{n+1} |I_{T'}|^{1/2}$$
for all sub-trees $T' \subset T \in \T$.
\end{definition}

The reader may easily verify that the modified energy is always dominated by the energy, and that we have the monotonicity property
$$
\modenergy_j((a_{P_j})_{\pv \in \Pv'} ) \leq
\modenergy_j((a_{P_j})_{\pv \in \Pv} )$$
whenever $\Pv' \subset \Pv$.  From duality we see that

\begin{lemma}\label{dual-energy}  Let the notation be as in Definition \ref{size-def}.  For any sequence $(a_{P_j})_{\pv \in \Pv}$, there exists a collection $\T$ of strongly $j$-disjoint trees, and complex co-efficients  $c_{P_j}$ for all $\pv \in \bigcup_{T \in \T} T$ such that
$$  \modenergy_j((a_{P_j})_{\pv \in \Pv} ) \sim  |\sum_{T \in \T} \sum_{\pv \in T} a_{P_j} \overline{c_{P_j}}|,$$
and such that
$$ \sum_{\pv \in T'} |c_{P_j}|^2 \lesssim \frac{|I_{T'}|}{\sum_{T \in \T} |I_T|}$$
for all $T \in \T$ and all sub-trees $T' \subseteq T$ of $T$.
\end{lemma}

\begin{proof}
Let $n$, $\T$ be an extremizer of \eqref{energy-mod-def}, and take $c_{P_j} := 2^{-n} (\sum_{T\in \T}|I_T|)^{-1/2} a_{P_j}$ for all 
$\pv \in \bigcup_{T \in \T} T$.
\end{proof}

The usual BMO norm can be written using an $L^2$ oscillation or an $L^1$ oscillation, and the two notions are equivalent thanks to the John-Nirenberg inequality.  The analogous statement for size is

\begin{lemma}\label{jn}
Let $\Pv$ be a finite collection of tri-tiles, $j=1,2,3$, and let $(a_{P_j})_{\pv \in \Pv}$ be a sequence of complex numbers. Then
\be{jn-est}
\size_j( (a_{P_j})_{\pv \in \Pv} ) \sim \sup_{T \subset \Pv}
\frac{1}{|I_T|} \| ( \sum_{\pv \in T} |a_{P_j}|^2 \frac{\chi_{I_\pv}}{|I_\pv|} )^{1/2} \|_{L^{1,\infty}(I_T)}
\end{equation}
where $T$ ranges over all trees in $\Pv$ which are $i$-trees for some $i \neq j$.
\end{lemma}

\begin{proof}
The  same as in \cite{mtt:walshbiest}, Lemma 4.2 .
\end{proof}

The following estimate is standard but we reproduce a proof in the Appendix for easy reference.  This is the main combinatorial tool needed to obtain estimates on \eqref{trilinear}.

\begin{proposition}\label{abstract}
Let $\Pv$ be a finite collection of tri-tiles, and for each $\pv \in \Pv$ and $j=1,2,3$ let $a^{(j)}_{P_j}$ be a complex number.  Then
\be{trilinear-est}
| \sum_{\pv \in \Pv} \frac{1}{|I_\pv|^{1/2}} a^{(1)}_{P_1} a^{(2)}_{P_2} a^{(3)}_{P_3}|
\lesssim \prod_{j=1}^3 
\size_j( (a^{(j)}_{P_j})_{\pv \in \Pv} )^{\theta_j}
\modenergy_j( (a^{(j)}_{P_j})_{\pv \in \Pv} )^{1-\theta_j}
\end{equation}
for any $0 \leq \theta_1, \theta_2, \theta_3 < 1$ with $\theta_1 + \theta_2 + \theta_3 = 1$, with the implicit constant depending on the $\theta_j$.
\end{proposition}

Note that this Proposition is stronger than that the corresponding statement (\cite{mtt:walshbiest}, Proposition 4.3) for the unmodified energy.  


Of course, in order to use Proposition \ref{abstract} we will need some estimates on size and energy.  In the rest of this section we give these estimates in the case when $a^{(j)}$ is given by \eqref{aj-def}.

We begin with a standard variant of Bessel's inequality for Fourier wave packets.

\begin{lemma}\label{l2}
Let $1\le j\le 3$, let $\T$ be a collection of strongly $j$-disjoint trees in $\Qv$, and for each $\qv \in T  \in \T$ let $c_{Q_j}$ be a complex number such that
\be{shibo}
\sum_{\qv \in \tilde T} |c_{Q_j}|^2 \lesssim A |I_{T'}|
\end{equation}
for all $T \in \T$ and sub-trees $T' \subseteq T$ of $T$, and some $A > 0$.  Then we have
$$ \| \sum_{T \in \T} \sum_{\qv \in T} c_{Q_j} \tilde \phi_{Q_j} \|_2 \lesssim (A \sum_{T \in \T} |I_T|)^{1/2}.$$
\end{lemma}

\begin{proof}
In the Walsh case this is immediate from Bessel's inequality. The argument in the Fourier case is more technical, however.

We may assume all trees in this lemma are sparse.

By squaring both sides, we reduce to showing that
$$ \sum_{T, T' \in \T} \sum_{\qv \in T} \sum_{\qv' \in T'}
|c_{Q_j}| |c_{Q'_j}| |\langle \tilde \phi_{Q_j}, \tilde \phi_{Q'_j} \rangle|
\lesssim A \sum_{T \in \T} |I_T|.$$
We may assume that $\omega_{Q_j} \cap \omega_{Q'_j} \neq \emptyset$ since the inner product vanishes otherwise.  By symmetry we may thus assume that $|\omega_{Q_j}|\le |\omega_{Q'_j}|$.

From the decay of the $\tilde \phi_{Q_j}$ we have
$$ |\langle \tilde \phi_{Q_j}, \tilde \phi_{Q'_j} \rangle|
\lesssim \frac{|I_{\qv'}|^{1/2}}{|I_\qv|^{1/2}}
(1 + \frac{\dist(I_\qv,I_{\qv'})}{|I_\qv|})^{-100}$$
so it suffices to show that
$$ \sum_{T, T' \in \T} \sum_{\qv \in T} \sum_{\qv' \in T': 
\omega_{Q_j} \cap \omega_{Q'_j}\neq \emptyset;
|\omega_{Q_j}| \le |\omega_{Q'_j}|}
|c_{Q_j}| |c_{Q'_j}| \frac{|I_{\qv'}|^{1/2}}{|I_\qv|^{1/2}}
(1 + \frac{\dist(I_\qv,I_{\qv'})}{|I_\qv|})^{-100}
\lesssim A \sum_{T \in \T} |I_T|.$$

Let us first consider the portion of the sum where $|\omega_{Q_j}| \sim |\omega_{Q'_j}|$.  
In this case we estimate $|c_{Q_j}| |c_{Q'_j}| \lesssim |c_{Q_j}|^2 + |c_{Q'_j}|^2$.  We treat the contribution of the first term $|c_{Q_j}|^2$, as the second is similar.  For each fixed $\qv$, 
there are only $O(1)$ many candidates $\omega $ to appear as $\omega_{Q'_j}$ 
satisfying all the above conditions, and for each fixed $\omega$
the $T', \qv'$ summations have disjoint spatial intervals $I_{\qv'}$.  One can then perform the $T', \qv'$ summations and estimate this contribution by
$$ \sum_{T \in \T} \sum_{\qv \in T} |c_{Q_j}|^2$$
which is acceptable by \eqref{shibo}.

It remains to consider the contribution when $|\omega_{Q_j}| \ll |\omega_{Q'_j}|$.  
From \eqref{shibo} applied to the singleton trees $\{\qv\}$, $\{\qv'\}$ we have
$$ |c_{Q_j}| \lesssim A^{1/2} |I_\qv|^{1/2}; \quad |c_{Q'_j}| \lesssim A^{1/2} |I_{\qv'}|^{1/2}.$$
It thus suffices to show that
$$ \sum_{T' \in \T} \sum_{\qv \in T} \sum_{\qv' \in T': \omega_{Q_j} \cap \omega_{Q'_j}\neq \emptyset; |I_\qv| \gg |I_{\qv'}|}
|I_{\qv'}| (1 + \frac{\dist(I_\qv,I_{\qv'})}{|I_\qv|})^{-100}
\lesssim |I_T|$$
for all trees $T \in \T$.

From the assumptions on $\qv$ and $\qv'$ and sparseness of the trees
we see that the tree $T'$ which contains $\qv'$ must be distinct from $T$.  By strong $j$-disjointness this implies that $I_{\qv'} \cap I_T = \emptyset$.  Also from strong $j$-disjointness we see that the $I_{\qv'}$ are disjoint.  We thus have
\bas
\sum_{\qv' \in T': \omega_{Q_j} \cap \omega_{Q'_j}\neq \emptyset; |I_\qv| > |I_{\qv'}|}
|I_{\qv'}| (1 + \frac{\dist(I_\qv,I_{\qv'})}{|I_\qv|})^{-100}
&\lesssim \int_{I_T^c} (1 + \frac{\dist(I_\qv,x)}{|I_\qv|})^{-100}\ dx\\
&\lesssim |I_\qv| (1 + \frac{\dist(I_\qv,I_T^c)}{|I_\qv|})^{-10}
\end{align*}
and the claim then follows by summing in $\qv$.
\end{proof}

As a consequence we have
 
\begin{lemma}\label{energy-lemma}
Let $j=1,2,3$, $f_j$ be a function in $L^2(\R)$, and let $\Pv$ be a finite collection of tri-tiles.  Then we have
\be{energy-lemma-est}
\modenergy_j((\langle f_j, \phi_{P_j} \rangle)_{P \in \P} ) \leq
\| f_j \|_2.
\end{equation}
\end{lemma}

\begin{proof}  
By Lemma \ref{dual-energy} we may find a collection of strongly $j$-disjoint trees $\T$, and complex co-efficients  $c_{P_j}$ for all $\pv \in \bigcup_{T \in \T} T$ such that
$$  \modenergy_j((\langle f_j, \phi_{P_j} \rangle)_{P \in \P} ) \sim
|\langle f_j, \sum_{T \in \T} \sum_{\pv \in T} c_{P_j} \phi_{P_j} \rangle|,$$
and such that
$$ \sum_{\pv \in T'} |c_{P_j}|^2 \lesssim \frac{|I_{T'}|}{\sum_{T \in \T} |I_T|}$$
for all $T \in \T$ and all sub-trees $T' \subseteq T$ of $T$.  The claim then follows from Cauchy-Schwarz and Lemma \ref{l2}.
\end{proof}

\begin{lemma}\label{size-lemma}
Let $j=1,2,3$, $E_j$ be a set of finite measure, $f_j$ be a function in $X(E_j)$, and let $\Pv$ be a finite collection of tri-tiles.  Then we have
\be{size-lemma-est}
\size_j( (\langle f_j, \phi_{P_j} \rangle)_{\pv \in \Pv} ) \lesssim
\sup_{\pv \in \Pv} \frac{\int_{E_j} \tilde \chi_{I_\pv}^{M}}{|I_\pv|}
\end{equation}
for all $M$, with the implicit constant depending on $M$.
\end{lemma}

\begin{proof}
This is essentially Lemma 7.8 in \cite{cct} (see also Lemma 4.5 of \cite{mtt:walshbiest}), but we give a proof here for completeness.

By Lemma \ref{jn} it suffices to show the estimate 
$$
\| (\sum_{\pv \in T} | \langle f_j, \phi_{P_j} \rangle |^2 \frac{\chi_{I_\pv}}{|I_\pv|})^{1/2}
\|_{L^{1,\infty}} \lesssim \int_{E_j} \tilde \chi_{I_T}^{M}
$$
for all $i \neq j$ and $i$-trees $T$.  

Fix $T$.  By frequency translation invariance we may assume that $\omega_{T,j}$ contains the origin.  

Let us first assume that $f_j$ is supported outside of $2I_T$.  From the decay of $\phi_{P_j}$ we have
\be{fp-decay}
|\langle f_j, \phi_{P_j} \rangle| \lesssim (\frac{|I_\pv|}{|I_T|})^M 
|I_T|^{-1/2} \int_{E_j} \tilde \chi_{I_T}^{M}
\end{equation}
Applying this estimate, we obtain
$$
\| (\sum_{\pv \in T} | \langle f_j, \phi_{P_j} \rangle |^2 \frac{\chi_{I_\pv}}{|I_\pv|})^{1/2}
\|_2 \lesssim 
|I_T|^{-1/2} \int_{E_j} \tilde \chi_{I_T}^{M}
$$
and the claim follows from H\"older.

Now suppose that $f_j$ is supported on $2I_T$.  It suffices to show that
\be{wt2} |\{ (\sum_{\pv \in T} | \langle f_j, \phi_{P_j} \rangle |^2 \frac{\tilde \chi_{I_\pv}^10}{|I_\pv|})^{1/2} 
\gtrsim \alpha \}| \lesssim \alpha^{-1} |E_j \cap 2I_T|.
\end{equation}
for all $\alpha > 0$.  

Since $\omega_{T,j}$ contains the origin, we see that the vector-valued operator
$$ f \mapsto ( \langle f, \phi_{P_j} \rangle \frac{\tilde \chi_{I_\pv}^10}{|I_\pv|} )_{\pv \in T}$$
is a Calder\'on-Zygmund operator, and is hence weak-type $(1,1)$; note that the $L^2$ boundedness of this operator follows from the almost orthogonality of the $\phi_{P_j}$.  (For more general $\omega_{T,j}$ this operator would be a modulated Calder\'on-Zygmund operator).  The claim follows from standard Calder\'on-Zygmund theory.
\end{proof}

In the next section we shall show how the above size and energy estimates can be combined with Proposition \ref{abstract} and the interpolation theory of the previous section to obtain Theorem \ref{bht}.  To prove the estimates for the trilinear operator $T_{\P,\Q}$ we need some more sophisticated size and energy estimates, which we will pursue after the proof of Theorem \ref{bht}.

\section{Proof of Theorem \ref{bht}}\label{bht-sec}

We now sketch a proof of Theorem \ref{bht}.  The proof here is standard, but we give it here for expository purposes, and also because we shall need Theorem \ref{bht} to prove the size and energy estimates needed for Theorem \ref{teorema1}.

We dualize $B$ into the trilinear form $\Lambda_{BHT}$ defined by
$$
\Lambda_{BHT}(f_1,f_2,f_3)(x):=\int_{\xi_1<\xi_2}
\widehat{f}_1(\xi_1)\widehat{f}_2(\xi_2) \widehat{f}_3(\xi_3)
\delta(\xi_1 + \xi_2 + \xi_3)\,d\xi_1 d\xi_2\ d\xi_3,$$
with $\delta$ denoting the Dirac delta. 

By standard discretization arguments as in Section \ref{discrete-sec} 
we may reduce the study of $\Lambda_{BHT}$ to that of 
discretized operators of the form 
\be{bht-disc-def}
\Lambda_{BHT,\Pv}
:= \sum_{\pv \in \Pv} \frac{1}{|I_\pv|^{1/2}}
\langle f_1, \phi_{P_1} \rangle
\langle f_2, \phi_{P_2} \rangle
\langle f_3, \phi_{P_3} \rangle
\end{equation}
where $\Pv$ is some finite collection of tri-tiles of rank 1.

We shall use the notation of Section \ref{interp-sec}, with the obvious modification for trilinear forms as opposed to quadrilinear forms.
From the interpolation theory in \cite{cct} it suffices to show that $\Lambda_{BHT,\Pv}$ is of restricted weak type $\alpha$ for all admissible 3-tuples $(\alpha_1, \alpha_2, \alpha_3)$ in the interior of the hexagon with vertices given by the six possible permutations of $(1,1/2,-1/2)$.  By symmetry and interpolation it suffices to prove restricted weak type $\alpha$ for admissible 3-tuples $\alpha$ arbitrarily close to $(1,1/2,-1/2)$, so that the bad index is 3.

Fix $\Pv$, $\alpha$ as above, and let $E_1$, $E_2$, $E_3$ be sets of finite measure.  We need to find a major subset $E'_3$ of $E_3$ such that
$$
|\Lambda_{BHT,\Pv}(f_1, f_2, f_3)|\lesssim |E|^{\alpha}
$$
for all functions $f_i\in X(E'_i)$, $i=1,2,3$.

Define the exceptional set $\Omega$ by
\[\Omega := \bigcup_{j=1}^{3}\{M\chi_{E_j}>C |E_j|/|E_3|\}\]
where $M$ is the dyadic Hardy-Littlewood maximal  function.
By the classical Hardy-Littlewood inequality, we have $|\Omega|<|E_3|/2$
if $C$ is a sufficiently large constant.  Thus if we set $E'_3 := E_3 \setminus \Omega$, then $E'_3$ is a major subset of $E_3$.

Let $f_i \in X(E'_i)$ for $i=1,2,3$.  We need to show
\be{choppy}
|\sum_{\pv \in \Pv} \frac{1}{|I_\pv|^{1/2}}
a^{(1)}_{P_1} a^{(2)}_{P_2} a^{(3)}_{P_3}| \lesssim |E|^\alpha
\end{equation}
where $a^{(j)}_{P_j}$ is defined by \eqref{aj-def}.

We shall make the assumption that
$$ 1 + \frac{\dist(I_\pv, \R \backslash \Omega)}{|I_\pv|} \sim 2^k$$
for all $\pv \in \Pv$, for some $k \geq 0$ independent of $\pv$, and prove \eqref{choppy} with an additional factor of $2^{-k}$ on the right-hand side.  If we can prove \eqref{choppy} in this special case with the indicated gain, then the general case of \eqref{choppy} follows by summing in $k$.

Fix $k$.  By the definition of $\Omega$ we have
$$ \frac{\int_{E_j} \tilde \chi_{I_\pv}^M}{|I_\pv|} \lesssim 2^k \frac{|E_j|}{|E_3|}$$
for all $\pv \in \Pv$ and $j=1,2$ and $M \gg 1$, while
$$ \frac{\int_{E'_3} \tilde \chi_{I_\pv}^M}{|I_\pv|} \lesssim 2^{(-M+C)k}$$
for all $\pv \in \Pv$ and $M \gg 1$ (of course, the implicit constant depends on $M$).

From Lemma \ref{size-lemma} we thus have
$$ \size_j( (a^{(j)}_{P_j})_{\pv \in \Pv} ) \lesssim 2^k \frac{|E_j|}{|E_3|}$$
for $j=1,2$, while
$$ \size_3( (a^{(3)}_{P_3})_{\pv \in \Pv} ) \lesssim 2^{-Mk}$$
for any $M$.

Also, from Lemma \ref{energy-lemma} and the fact that $f_j \in X(E'_j)$ we have
$$ \modenergy_j( (a^{(j)}_{P_j})_{\pv \in \Pv} ) \lesssim |E_j|^{1/2}$$
for $j=1,2,3$.

From Proposition \ref{abstract} and a suitably large choice of $M$ 
depending on the $\theta_j$ we thus have
$$
|\sum_{\pv \in \Pv} \frac{1}{|I_\pv|^{1/2}}
a^{(1)}_{P_1} a^{(2)}_{P_2} a^{(3)}_{P_3}| \lesssim 2^{-k}
\prod_{j=1}^3 |E_j|^{(1-\theta_j)/2} (\frac{|E_j|}{|E_3|})^{\theta_j} $$
for any $0 < \theta_1, \theta_2, \theta_3 < 1$ such that $\theta_1 + \theta_2 + \theta_3 = 1$.  The claim then follows by choosing $\theta_1 := 2\alpha_1-1$, $\theta_2 := 2\alpha_2-1$, and $\theta_3 := 2\alpha_3 + 1$; note that there exist choices of $\alpha$ arbitrarily close to $(1,1/2,-1/2)$ for which the constraints on $\theta_1, \theta_2, \theta_3$ are satisfied.  
This concludes the proof of Theorem \ref{bht}.

\section{Energy estimates}\label{energy-sec}

The purpose of this section is to prove some additional energy estimates in the spirit of Lemma \ref{energy-lemma}.

From Lemma \ref{l2} and the Cauchy-Schwarz inequality we have

\begin{corollary}\label{corollary}
Let $\T$ be a collection of strongly $3$-disjoint trees in $\Qv$, 
and for each $\qv \in T  \in \T$ let $c_{Q_3}$ be a complex number such that
\be{ask-1}
\sum_{\qv \in \tilde T} |c_{Q_3}|^2 \lesssim \frac{|I_{\tilde T}|}{\sum_{T \in \T} |I_T|}
\end{equation}
for all $\tilde T \subseteq T \in \T$.  Also, let $\T'$ be a collection of $1$-disjoint trees in $\Pv$, and for each $\pv \in T' \in \T'$ let $d_{P_1}$ be a complex number such that
\be{ask-2}
\sum_{\pv \in \tilde T} |d_{P_1}|^2 \lesssim \frac{|I_{\tilde T}|}{\sum_{T' \in \T'} |I_{T'}|}
\end{equation}
for all $\tilde T \subseteq T' \in \T'$.  Then we have
\be{job-0}
|\sum_{T' \in \T'} \sum_{\pv \in T'}
\sum_{T \in \T} \sum_{\qv \in T}
c_{Q_3} d_{P_1} \langle \phi_{P_1}, \tilde \phi_{Q_3} \rangle| \lesssim 1.
\end{equation}
\end{corollary}

In the rest of this section we shall prove \eqref{job-0} when the constraint $\omega_{Q_3} \subseteq \omega_{P_1}$ has been inserted into the summation.  In \cite{mtt:walshbiest} this was accomplished by a geometric lemma (\cite{mtt:walshbiest}, Lemma 6.1) which allowed one to decouple the $\omega_{Q_3} \subseteq \omega_{P_1}$ constraint assuming a priori that $P_1 \cap Q_3 \neq \emptyset$.  In the Walsh case this assumption was reasonable, however in the Fourier case we do not have perfect orthogonality in space and so we can only assume that $\omega_{P_1} \cap \omega_{Q_3} \neq \emptyset$ a priori.  However, we still have the following weaker analogue of \cite{mtt:walshbiest}, Lemma 6.1 when the $\qv$ tri-tiles are constrained to a tree.

\begin{lemma}\label{biest-trick}
Let $T \subset \Qv$ be a (sparse) $i$-tree for some $i=1,2$, and
define the collection $\Pv' \subset \Pv$ of tri-tiles by
$$ \Pv' := \{ \pv \in \Pv: \omega_{Q_3} \subseteq \omega_{P_1} \hbox{ for some } \qv \in T\}.$$
Then, if $\qv \in T$ and $\pv \in \Pv$ are such that $\langle \phi_{P_1}, \tilde \phi_{Q_3} \rangle \neq 0$, then
$$ \omega_{Q_3} \subseteq \omega_{P_1} \iff \pv \in \Pv'.$$
\end{lemma}

\begin{proof}
If $\langle \phi_{P_1}, \tilde \phi_{Q_3} \rangle \neq 0$, then $\omega_{Q_3} \cap \omega_{P_1} \neq \emptyset$.  The claim then follows from the sparseness of $T$. 
\end{proof}

\begin{lemma}\label{split}  Let the notation be as in Corollary \ref{corollary}. Then we have
\be{job}
|\sum_{T' \in \T'} \sum_{\pv \in T'}
\sum_{T \in \T} \sum_{\qv \in T: \omega_{Q_3} \subseteq \omega_{P_1}}
c_{Q_3} d_{P_1} \langle \phi_{P_1}, \tilde \phi_{Q_3} \rangle| \lesssim 1.
\end{equation}
\end{lemma}

\begin{proof}

We can divide into the cases
\be{it-itp}
\sum_{T \in \T} |I_T| \lesssim \sum_{T' \in \T'} |I_{T'}|
\end{equation}
and
\be{itp-it}
\sum_{T' \in \T'} |I_{T'}| \lesssim \sum_{T \in \T} |I_{T}|.
\end{equation}
In the case \eqref{itp-it} we use \eqref{job-0} to reduce \eqref{job} to
$$
|\sum_{T \in \T} \sum_{\qv \in T}
\sum_{T' \in \T'} \sum_{\pv \in T': \omega_{P_1} \subsetneq \omega_{Q_3}}
c_{Q_3}  d_{P_1} 
\langle 
\phi_{P_1},
\tilde \phi_{Q_3} 
\rangle| \lesssim 1.
$$
But the proof of this estimate is essentially the same as \eqref{job} with the roles of $P$ and $Q$ reversed.
Thus it suffices to prove \eqref{job} under the assumption \eqref{it-itp}.

We first consider the set of all pairs $\omega_{Q_3}\subset\omega_{P_1}$
such that $\omega_{T,3}\not\subset 2\omega_{P_1}$, where $T$ is the
tree in $\T$ containing $Q$. These constraints imply
$10^9|\omega_{Q_3}|>|\omega_{P_1}|$. By splitting into $O(1)$
cases we may assume that the ratio between $|\omega_{Q_3}|$ and 
$|\omega_{P_1}|$
is fixed. We may also assume that the distance of $I_{Q}$ and $I_P$
is $2^k|I_Q|$ for some fixed $k$, provided we prove the final estimate
with an extra factor of $2^{-k}$. However, then we have
\[\langle \phi_{P_1}, \tilde \phi_{Q_3} \rangle\neq 0\]
only for a bounded number of essentially unique $P$ for any given $Q$, 
and for those $P$ we have 
\[\langle \phi_{P_1}, \tilde \phi_{Q_3} \rangle \lesssim 2^{-k}\]
Hence we can estimate the corresponding piece of 
\eqref{job} using Cauchy-Schwarz by  
\[ 2^{-k} ( \sum_{T \in \T} \sum_{\qv \in T}|c_{Q_3}|^2 )^{\frac 12}
(\sum_{T' \in \T'} \sum_{\pv \in T'}|d_{P_1}|^2)^{\frac 12}
\lesssim 2^{-k}\]

Now we consider the pairs $(P,Q)$ with 
$\omega_{T,3}\subset 2\omega_{P_1}
$, where $T$ is the tree containing $Q$.
We estimate the corresponding part of the 
left-hand side of \eqref{job} by
$$
\sum_{T \in \T} \sum_{\qv \in T} |c_{Q_3}|
|\langle \sum_{T' \in \T'} \sum_{\pv \in T': \omega_{Q_3} \subseteq \omega_{P_1},\omega_{T,3}\subset 2\omega_{P_1}
} 
d_{P_1} \phi_{P_1}, \tilde \phi_{Q_3} \rangle|.
$$
By \eqref{it-itp}, it suffices to show that
\[
\sum_{\qv \in T} |c_{Q_3}|
|\langle \sum_{T' \in \T'} \sum_{\pv \in T': \omega_{Q_3} \subseteq \omega_{P_1},\omega_{T,3}\subset 2\omega_{P_1}} d_{P_1} \phi_{P_1}, \tilde \phi_{Q_3} \rangle|\]
\be{job-1}
\lesssim \frac{|I_T|}{(\sum_{T \in \T} |I_T|)^{1/2}}
\frac{1}{(\sum_{T' \in \T'} |I_{T'}|)^{1/2}}
\end{equation}
for each $T \in \T$.

Fix $T$. Let us first estimate the contribution of the case when $I_\pv \cap 2I_T \neq \emptyset$.  Define the collection $\Pv_T$ by
\be{pvt-def}
\Pv_T := \{ \pv \in \bigcup_{T' \in \T'} T': I_\pv \cap 2I_T \neq \emptyset; 
\quad \omega_{T,3}\subset 2\omega_{P_1};
\quad \omega_{Q_3} \subseteq \omega_{P_1} \hbox{ for some } Q_3 \in T \}.
\end{equation}
By Lemma \ref{biest-trick} we may rewrite the contribution of this case to \eqref{job-1} as
$$\sum_{\qv \in T} |c_{Q_3}|
|\langle h_T, \tilde \phi_{Q_3} \rangle|.$$
where
$$ h_T := \sum_{\pv \in \Pv_T} d_{P_1} \phi_{P_1}.$$
By Cauchy-Schwarz and \eqref{ask-1} we can bound the previous by
$$
 (\frac{|I_T|}{\sum_{T \in \T} |I_T|})^{1/2}
(\sum_{\qv \in T} |\langle h_T, \tilde \phi_{Q_3} \rangle|^2)^{1/2}.$$
The $\tilde \phi_{Q_3}$ are almost orthogonal as $\qv \in T$ varies, so we can bound this by
$$
 (\frac{|I_T|}{\sum_{T \in \T} |I_T|})^{1/2}
\| h_T \|_2.$$
It thus suffices to prove 
\be{ht-est}
\| h_T\|_2 \lesssim \frac{|I_T|^{1/2}}{(\sum_{T' \in \T'} |I_{T'}|)^{1/2}}.
\end{equation}
We write the left-hand side as
$$ \| \sum_{T' \in \T'} \sum_{\pv \in T' \cap \Pv_T} d_{P_1} \phi_{P_1} \|_2.$$
We now consider each $\{\pv\}$ with $\pv\in\Pv_T$ as a tree by itself.
By strongly 1- disjointness of the tree $T'$ we see that
$I_P$ with $\pv\in \Pv_T$ are pairwise disjoint. Moreover, they are contained
in $3I_T$.  
In particular we have $\sum_{\pv\in \Pv_T} |I_{P}| \lesssim |I_T|$, 
and the claim \eqref{ht-est} follows from Lemma \ref{l2}.  This concludes the treatment of the case $ I_\pv \cap 2I_T \neq \emptyset$.

To finish the estimation of \eqref{job-1} it remains to treat the contribution of the case $I_\pv \cap 2^kI_T = \emptyset$ and $I_\pv \cap 2^{k+1}I_T \neq \emptyset$ for each $k > 0$, with an additional factor of $2^{-k}$ on the right hand side.

Fix $k > 0$.  In this case we use the crude estimate 
$$|c_{Q_3}| \lesssim \frac{|I_\qv|^{1/2}}{(\sum_{T \in \T} |I_T|)^{1/2}}$$
from \eqref{ask-1}, and reduce to showing that
$$
|\langle \sum_{T' \in \T'} \sum_{\pv \in T': \omega_{Q_3} \subseteq \omega_{P_1}; \omega_{T,3}\subset 2\omega_{P_1}; I_\pv \cap 2^k I_T = \emptyset; I_\pv \cap 2^{k+1} I_T \neq \emptyset} 
d_{P_1} \phi_{P_1}, \tilde \phi_{Q_3} \rangle|$$
$$\lesssim 2^{-k} (\frac{|I_\qv|}{|I_T|})^{10}
\frac{|I_T|^{1/2}}{(\sum_{T' \in \T'} |I_{T'}|)^{1/2}}$$
for each $\qv \in \T$.

Fix $Q$.  We split $\tilde \phi_{Q_3}$ into $\tilde \phi_{Q_3} \chi_{2^{k-1} I_T}$ and
$\tilde \phi_{Q_3} (1 - \chi_{2^{k-1} I_T})$.  To control the former contribution we use the crude estimates
$$ |d_{P_1}| \lesssim \frac{|I_\pv|^{1/2}}{(\sum_{T' \in \T'} |I_{T'}|)^{1/2}}$$
from \eqref{ask-2} and
$$ |\langle  \phi_{P_1}, \tilde \phi_{Q_3} \chi_{2^{k-1} I_T} \rangle| \lesssim 
2^{-100k} (\frac{|I_\qv|}{|I_T|})^{100} (\frac{|I_\pv|}{|I_T|})^{100},$$
and sum crudely in $\pv$.  To control the latter contribution we observe that 
$\tilde \phi_{Q_3} \chi_{2^{k-1} I_T}$ has an $L^2$ norm of $O(2^{-100k} (\frac{|I_\qv|}{|I_T|})^{100})$, so it suffices to show that 
$$ \| \sum_{T' \in \T'} \sum_{\pv \in T': \omega_{Q_3} \subseteq \omega_{P_1}; I_\pv \cap 2^k I_T = \emptyset; 
I_\pv \cap 2^{k+1} I_T \neq \emptyset;
\omega_{T,3}\subset 2\omega_{P_1}} 
d_{P_1} \phi_{P_1} \|_2 \lesssim 2^{50k} (\frac{|I_T|}{|I_\qv|})^{50} \frac{|I_T|^{1/2}}{(\sum_{T' \in \T'} |I_{T'}|)^{1/2}}.$$
But this follows by repeating the proof of \eqref{ht-est}.
\end{proof}

From Lemma \ref{split} and Lemma \ref{dual-energy} we have

\begin{corollary}\label{split-cor}
Let $\T$ be a collection of $3$-disjoint trees in $\Qv$, and for each $\qv \in T  \in \T$ let $c_{Q_3}$ be a complex number such that \eqref{ask-1} holds for all $\tilde T \subseteq T \in \T$.  Then 
$$
\modenergy_1(
(\sum_{T \in \T} \sum_{\qv \in T: \omega_{Q_3} \subseteq \omega_{P_1}}
c_{Q_3} \langle \phi_{P_1}, \tilde \phi_{Q_3} \rangle)_{\pv \in \Pv} ) \lesssim 1.
$$
\end{corollary}

\section{Additional size and energy estimates}

In the expression $\Lambda_{\Pv,\Qv}$ the $Q$ tiles in the inner summation have a narrower frequency interval, and hence a wider spatial interval, than the $P$ tiles in the outer summation.  Thus the inner summation has a poorer spatial localization than the outer sum.  It shall be convenient to reverse the order of summation so that the inner summation is instead more strongly localized spatially than the outer summation.  Specifically, we rewrite $\Lambda_{\Pv,\Qv}$ as 
$$
\Lambda_{\Pv,\Qv}(f_1,f_2,f_3,f_4)=
\sum_{\qv\in\Qv}\frac{1}{|I_\qv|^{1/2}} a^{(1)}_{Q_1} a^{(2)}_{Q_2} a^{(3)}_{Q_3}
$$
where
\begin{equation}\label{rightform}
\begin{split}
a^{(1)}_{Q_1} &:= \langle f_1,\tilde \phi_{Q_1} \rangle\\
a^{(2)}_{Q_2} &:= \langle f_2,\tilde \phi_{Q_2} \rangle\\
a^{(3)}_{Q_3} &:= \sum_{\pv\in\Pv\,;\,\omega_{Q_3}\subseteq \omega_{P_1}}
\frac{1}{|I_\pv|^{1/2}}
\langle f_3,\phi_{P_2} \rangle
\langle f_4,\phi_{P_3} \rangle
\langle \phi_{P_1}, \tilde \phi_{Q_3} \rangle. 
\end{split}
\end{equation}

The purpose of this section is to prove 
analogues of Lemma \ref{size-lemma} and 
Lemma \ref{energy-lemma} for $a^{(3)}_{Q_3}$.

\begin{lemma}\label{bht-size}
Let $E_j$ be sets of finite measure and $f_j$ be functions in $X(E_j)$ for $j=3,4$.  Then we have
\be{size-bht-est}
\size_3((a^{(3)}_{Q_3})_{\qv\in \Qv}) \lesssim
\sup_{\qv \in \Qv} \left( \frac{\int_{E_3} \tilde \chi_{I_\qv}^M}{|I_\qv|}\right)^{1-\theta} \left( \frac{\int_{E_4} \tilde \chi_{I_\qv}^M}{|I_\qv|}\right)^\theta
\end{equation}
for any $0 < \theta < 1$ and $M > 0$, with the implicit constant depending on $\theta, M$.
\end{lemma}

\begin{proof}
By Lemma \ref{jn} it suffices to show that
$$
\| (\sum_{\qv \in T} |a^{(3)}_{Q_3}|^2 \frac{\chi_{I_\qv}}{|I_\qv|})^{1/2} \|_{L^{1,\infty}(I_T)}
\lesssim |I_T|
\sup_{\qv \in \Qv} (\frac{\int_{E_3} \tilde \chi_{I_\qv}^M}{|I_\qv|})^{1-\theta} (\frac{\int_{E_4} \tilde \chi_{I_\qv}^M}{|I_\qv|})^\theta
$$
for any $i = 1,2$ and any $i$-tree $T$.  We may assume (as in the proof of Lemma \ref{size-lemma}) that $T$ contains its top $P_T$, in which case we may reduce to
\be{weak-sap}
\| (\sum_{\qv \in T} |a^{(3)}_{Q_3}|^2 \frac{\chi_{I_\qv}}{|I_\qv|})^{1/2} \|_{L^{1,\infty}(I_T)}
\lesssim 
(\int_{E_3} \tilde \chi_{I_T}^M)^{1-\theta} (\int_{E_4} \tilde \chi_{I_T}^M)^\theta.
\end{equation}

Fix $T$.  To prove \eqref{weak-sap}, first consider the relatively easy case when $f_3$ vanishes on $5I_T$.  In this case we shall prove the stronger estimate 
\begin{align}
|a^{(3)}_{Q_3}| &\lesssim 
|I_\qv|^{-1/2}
(\int_{E_3} \tilde \chi_{I_\qv}^M)^{1-\theta} (\int_{E_4} \tilde \chi_{I_\qv}^M)^\theta\label{weak-targ}\\
&\lesssim
|I_\qv|^{-1/2} (\frac{|I_\qv|}{|I_T|})^{M(1-\theta)}
(\int_{E_3} \tilde \chi_{I_T}^M)^{1-\theta} (\int_{E_4} \tilde \chi_{I_T}^M)^\theta\nonumber
\end{align}
for all $\qv \in T$; the claim \eqref{weak-sap} then follows by square-summing in $\qv$.  

We now prove \eqref{weak-targ}.  Fix $\qv \in \T$.  By \eqref{rightform} and \eqref{decay} we may estimate
$$
|a^{(3)}_{Q_3}|
\lesssim |I_\qv|^{-1/2}
\sum_{\pv\in\Pv\,;\,\omega_{Q_3}\subseteq \omega_{P_1}}
|\langle f_3,\phi_{P_2} \rangle|
|\langle f_4,\phi_{P_3} \rangle|
\int
\frac{\tilde \chi_{I_\pv}^{100M}}{|I_\pv|} \tilde \chi_{I_\qv}^{100M}.
$$
Interchanging the sum and integral and applying Cauchy-Schwarz we thus have
$$
|a^{(3)}_{Q_3}|
\lesssim 
\int |S_2 f_3| |S_3 f_4| \tilde \chi_{I_\qv}^{100M}
$$
where for $j=2,3$, the square function $S_j$ is the vector-valued quantity
$$  S_j f := 
(\langle f,\phi_{P_j} \rangle
\frac{\tilde \chi_{I_\pv}^{50M}}{|I_\pv|^{\frac 12}})_{\pv\in\Pv\,;\,\omega_{Q_3}\subseteq \omega_{P_1}}.$$
To show \eqref{weak-targ}, it thus suffices by H\"older to prove the weighted square-function estimate 
\be{sj-prof}
\| S_j f \|_{L^p(\tilde \chi_{I_\qv}^{100M}\ dx)}
\lesssim \| f \|_{L^p(\tilde \chi_{I_\qv}^M\ dx)}
\end{equation}
for all $1 < p < \infty$ and $j=2,3$.   But this follows since $S_j$ is a modulated Calder\'on-Zygmund operator\footnote{This is ultimately because the constraint $\omega_{Q_3} \subseteq \omega_{P_1}$, combined with the rank 1 hypothesis, forces the intervals $\omega_{P_j}$ to be lacunary around some frequency $\xi$.} whose kernel $K_j(x,y)$ decays like
$O(|I_\qv|^{-1} (|x-y|/|I_\qv|)^{-50M})$ for all $|x-y| \gg |I_\qv|$.  

This proves \eqref{weak-sap} when $f_3$ vanishes on $5I_T$.  A similar argument gives \eqref{weak-sap} when $f_4$ vanishes on $5I_T$.  We may thus reduce to the case when $f_3$, $f_4$ are both supported on $5I_T$.  We may then assume that $E_3, E_4 \subset 5I_T$.

Define the collection $\Pv' \subset \Pv$ of tri-tiles by
$$ \Pv' := \{ \pv \in \Pv: \omega_{Q_3} \subseteq \omega_{P_1} \hbox{ for some } \qv \in T\}.$$
From Lemma \ref{biest-trick} and \eqref{rightform} we have
\bas
a^{(3)}_{Q_3} &= \sum_{\pv\in\Pv'} \frac{1}{|I_\pv|^{1/2}}
\langle f_3,\phi_{P_2} \rangle
\langle f_4,\phi_{P_3} \rangle
\langle \phi_{P_1}, \tilde \phi_{Q_3} \rangle\\
&= \langle B_{\Pv'}(f_3,f_4), \tilde \phi_{Q_3} \rangle
\end{align*}
where
$$ B_{\Pv'}(f_3,f_4) :=
\sum_{\pv\in\Pv'} \frac{1}{|I_\pv|^{1/2}}
\langle f_3,\phi_{P_2} \rangle
\langle f_4,\phi_{P_3} \rangle \phi_{P_1}.$$

To prove \eqref{weak-sap} it thus suffices to show that
$$
\| (\sum_{\qv \in T} |\langle B_{\Pv'}(f_3,f_4), \phi_{Q_3} \rangle|^2 \frac{\tilde \chi_{I_\qv}^{100}}{|I_\qv|})^{1/2} \|_{L^{1,\infty}}
\lesssim 
|E_3|^{1-\theta} |E_4|^\theta.
$$
The vector-valued operator 
\be{vec}
f \mapsto (\langle f, \phi_{Q_3} \rangle \frac{\tilde \chi_{I_\qv}^{50}}{|I_\qv|^{\frac 12}})_{\qv \in T}
\end{equation}
is a modulated Calder\'on-Zygmund operator, so it suffices to show that
\be{jolt}
\| B_{\Pv'}(f_3,f_4) \|_1 \lesssim \| f_3 \|_{1/(1-\theta)} \|f_4 \|_{1/\theta}.
\end{equation}
But this follows from Theorem \ref{bht} (or more precisely, the analogue of Theorem \ref{bht} for the discretized operator $B_{\Pv'}$).  This finishes the proof of \eqref{size-bht-est}.
\end{proof}

The analogue of Lemma \ref{energy-lemma} is

\begin{lemma}\label{bht-energy}
Let $E_j$ be sets of finite measure and $f_j$ be functions in $X(E_j)$ for $j=3,4$.  Then we have
\be{energy-bht-est-2}
\modenergy_3((a^{(3)}_{Q_3})_{\qv \in \Qv}) \lesssim
\left(|E_4|^{1/2} \sup_{\pv \in \Pv} \frac{\int_{E_3} \tilde \chi_{I_\pv}^M}{|I_\pv|}\right)^{1-\theta}
\left(|E_3|^{1/2} \sup_{P \in \P} \frac{\int_{E_4} \tilde \chi_{I_\pv}^M}{|I_\pv|}\right)^{\theta}
\end{equation}
for any $0 < \theta < 1$ and $M>0$, with the implicit constant depending on $\theta$, $M$.  In particular, we have
\be{energy-bht-est}
\modenergy_3((a^{(3)}_{Q_3})_{\qv \in \Qv}) \lesssim
|E_4|^{(1-\theta)/2} |E_3|^{\theta/2}
\end{equation}
for any $0 < \theta < 1$, with the implicit constant depending on $\theta$.
\end{lemma}

\begin{proof}
By Lemma \ref{dual-energy}, it suffices to show that
\be{energy-1}
|\sum_{T \in \T} \sum_{\qv \in T} a_{Q_3} \overline{c_{Q_3}}| \lesssim
\left(|E_4|^{1/2} \sup_{\pv \in \Pv} \frac{\int_{E_3} \tilde \chi_{I_\pv}^M}{|I_\pv|}\right)^{1-\theta}
\left(|E_3|^{1/2} \sup_{P \in \P} \frac{\int_{E_4} \tilde \chi_{I_\pv}^M}{|I_\pv|}\right)^{\theta}
\end{equation}
for all collections $\T$ of strongly $3$-disjoint trees, and all co-efficients $c_{Q_3}$ such that
$$ \sum_{\qv \in \tilde T} |c_{Q_3}|^2 \sim \frac{|I_{\tilde T}|}{\sum_{T \in \T} |I_T|}$$
for all $\tilde T \subseteq T \in \T$.

Fix $\T$, $c_{Q_3}$.  By \eqref{rightform} we may write the left-hand side of \eqref{energy-1} as 
$$
|\sum_{\pv\in\Pv}
\frac{1}{|I_\pv|^{1/2}}
b^{(1)}_{P_1}
\langle f_3,\phi_{P_2} \rangle
\langle f_4,\phi_{P_3} \rangle|$$
where
\be{b-def}
b^{(1)}_{P_1}
:= 
\sum_{T \in \T} \sum_{\qv \in T: \omega_{Q_3} \subseteq \omega_{P_1}}
\langle \phi_{P_1}, c_{Q_3} \tilde \phi_{Q_3} \rangle. 
\end{equation}

The claim then follows from Proposition \ref{abstract}, Corollary \ref{split-cor}, Lemma \ref{energy-lemma}, and Lemma \ref{size-lemma}.
\end{proof}

\section{Proof of Theorem \ref{teorema1} for $A_5, \ldots, A_{12}$}

Let $\alpha=(\alpha_1,\alpha_2,\alpha_3,\alpha_4)$ admissible tuples near $A_i$ for some $5 \leq i \leq 12$.  We will only consider those vertices with bad index 1 (i.e. $A_9, \ldots, A_{12}$) as the other four vertices can be done similarly.  Thus $\alpha$ has bad index $1$. Let us also
fix arbitrary sets $E_1,E_2,E_3,E_4$ of finite measure. 

As in the proof of Theorem \ref{bht}, we define 
\[\Omega := \bigcup_{j=1}^{4}\{M\chi_{E_j}>C|E_j|/|E_1|\}\]
for a large constant $C$, and set $E'_1 := E_1 \setminus \Omega$.  We now fix $f_i \in X(E'_i)$ for $i=1,2,3,4$.  Our task is then to show

\begin{equation}\label{weakin}
|\sum_{\qv\in\Qv}\frac{1}{|I_\qv|^{1/2}} a^{(1)}_{Q_1} a^{(2)}_{Q_2} a^{(3)}_{Q_3}
|\lesssim |E|^\alpha
\end{equation}
where the $a^{(j)}_{Q_j}$ are defined by \eqref{rightform}.

As before, we may make the assumption that
$$ 1 + \frac{\dist(I_\qv, \R \backslash \Omega)}{|I_\qv|} \sim 2^k$$
for all $\qv \in \Qv$ and for some $k \geq 0$ independent of $\qv$, provided that we gain a factor such as $2^{-k}$ on the right-hand side of \eqref{weakin}.
As before, we then have
$$ \frac{\int_{E_j} \tilde \chi_{I_\qv}^M}{|I_\qv|} \lesssim 2^k \frac{|E_j|}{|E_1|}$$
for all $\qv \in \Qv$ and $j=2,3,4$ and $M \gg 1$, while
$$ \frac{\int_{E'_1} \tilde \chi_{I_\qv}^M}{|I_\qv|} \lesssim 2^{(-M+C)k}$$
for all $\qv \in \Qv$ and $M \gg 1$.  

With this assumption we have
From Lemma \ref{size-lemma} and Lemma \ref{bht-size} we thus have
\bas
\size_1( (a^{(1)}_{Q_1})_{\qv \in \Qv} ) &\lesssim 2^{-(M+C)k}\\
\size_2( (a^{(2)}_{Q_2})_{\qv \in \Qv} ) &\lesssim 2^k \frac{|E_2|}{|E_1|}\\
\size_3( (a^{(3)}_{Q_3})_{\qv \in \Qv} ) &\lesssim 2^k \frac{|E_3|^{1-\theta} |E_4|^\theta}{|E_1|}
\end{align*}
for some $0 < \theta < 1$ which we will choose later.  Similarly, from Lemma \ref{energy-lemma}, \eqref{energy-bht-est} and the hypotheses $f_j \in X(E_j)$ we have
\bas
\modenergy_1( (a^{(1)}_{Q_1})_{\qv \in \Qv} ) &\lesssim |E_1|^{1/2}\\
\modenergy_2( (a^{(2)}_{Q_2})_{\qv \in \Qv} ) &\lesssim |E_2|^{1/2}\\
\modenergy_3( (a^{(3)}_{Q_3})_{\qv \in \Qv} ) &\lesssim |E_3|^{(1-\theta)/2} |E_4|^{\theta/2}.
\end{align*}

By Proposition \ref{abstract} and choosing $M$ sufficiently large we can thus bound the left-hand side of \eqref{weakin} by 
$$
2^{-k}
\frac{|E_1|^{(1+\theta_1)/2} |E_2|^{(1+\theta_2)/2} (|E_3|^{1-\theta} |E_4|^\theta)^{(1+\theta_3)/2} }{|E_1|}$$
for $0 < \theta_1, \theta_2, \theta_3 < 1$ such that $\theta_1 + \theta_2 + \theta_3 = 1$.  The claim follows by setting $\theta_1 := 2\alpha_1 + 1$, $\theta_2 := 2\alpha_2 - 1$, $\theta_3 := 2(\alpha_3 + \alpha_4) - 1$, and $\theta := \alpha_4 / (\alpha_3 + \alpha_4)$; the reader may verify that the constraints on $\theta_1, \theta_2, \theta_3, \theta$ can be obeyed for $\alpha$ arbitrarily close to $A_9, A_{10}, A_{11}, A_{12}$.

\section{Proof of Theorem \ref{teorema1} for $A_1, A_2, A_3, A_4$}

Let $\alpha=(\alpha_1,\alpha_2,\alpha_3,\alpha_4)$ admissible tuples near $A_i$ for some $1 \leq i \leq 4$.  We will only consider those vertices with bad index 4 (i.e. $A_1, A_{2}$) as the other two vertices can be done similarly.  Thus $\alpha$ has bad index $4$. Let us also
fix arbitrary sets $E_1,E_2,E_3,E_4$ of finite measure. 

As before, we define 
\[\Omega := \bigcup_{j=1}^{4}\{M\chi_{E_j}>C|E_j|/|E_4|\}\]
for a large constant $C$, and set $E'_4 := E_4 \setminus \Omega$.  We now fix $f_i \in X(E'_i)$ for $i=1,2,3,4$.  Our task is then to show

\begin{equation}\label{weakin-2}
|\sum_{\qv\in\Qv}\frac{1}{|I_\qv|^{1/2}} a^{(1)}_{Q_1} a^{(2)}_{Q_2} a^{(3)}_{Q_3}
|\lesssim |E|^\alpha
\end{equation}
where the $a^{(j)}_{Q_j}$ are defined by \eqref{rightform}.

Recall that $a^{(3)}_{Q_3}$ is defined by
$$
a^{(3)}_{Q_3} := \sum_{\pv\in\Pv\,;\,\omega_{Q_3}\subseteq \omega_{P_1}}
\frac{1}{|I_\pv|^{1/2}}
\langle f_3,\phi_{P_2} \rangle
\langle f_4,\phi_{P_3} \rangle
\langle \phi_{P_1}, \phi_{Q_3} \rangle. 
$$

We may make the assumptions that
$$ 1 + \frac{\dist(I_\qv, \R \backslash \Omega)}{|I_\qv|} \sim 2^k$$
for all $\qv \in \Qv$ and for some $k \geq 0$ independent of $\qv$, and that
$$ 1 + \frac{\dist(I_\pv, \R \backslash \Omega)}{|I_\pv|} \sim 2^{k'}$$
for all $\pv \in \Pv$ and for some $k' \geq 0$ independent of $\pv$, 
provided that we gain a factor such as $2^{-k-k'}$ on the right-hand side of \eqref{weakin}.

As before, we then have
$$ \frac{\int_{E_j} \tilde \chi_{I_\qv}^M}{|I_\qv|} \lesssim 2^k \frac{|E_j|}{|E_4|}$$
for all $\qv \in \Qv$ and $j=1,2,3$ and $M \gg 1$, while
$$ \frac{\int_{E'_4} \tilde \chi_{I_\qv}^M}{|I_\qv|} \lesssim 2^{(-M+C)k}$$
for all $\qv \in \Qv$ and $M \gg 1$.  From Lemma \ref{size-lemma} we thus have
\bas
\size_1( (a^{(1)}_{Q_1})_{\qv \in \Qv} ) &\lesssim 2^k \frac{|E_1|}{|E_4|}\\
\size_2( (a^{(2)}_{Q_2})_{\qv \in \Qv} ) &\lesssim 2^k \frac{|E_2|}{|E_4|}.
\end{align*}
From Lemma \ref{bht-size} and the crude estimate $\int_{E_j} \tilde \chi_{I_\pv}^M \leq |I_\pv|$ we also have
$$
\size_3( (a^{(3)}_{Q_3})_{\qv \in \Qv} ) \lesssim 2^{-(M+C)k}.$$

From Lemma \ref{energy-lemma} we have
\bas
\modenergy_1( (a^{(1)}_{Q_1})_{\qv \in \Qv} ) &\lesssim |E_1|^{1/2}\\
\modenergy_2( (a^{(2)}_{Q_2})_{\qv \in \Qv} ) &\lesssim |E_2|^{1/2}.
\end{align*}

Finally, from the definition of $k'$ we have
$$ \frac{\int_{E_3} \tilde \chi_{I_\pv}^M}{|I_\pv|} \lesssim 2^k \frac{|E_3|}{|E_4|}$$
$$ \frac{\int_{E'_4} \tilde \chi_{I_\pv}^M}{|I_\pv|} \lesssim 2^{(-M+C)k}$$
for all $\pv \in \Pv$.  By \eqref{energy-bht-est-2} we thus have
\bas
\modenergy_3( (a^{(3)}_{Q_3})_{\qv \in \Qv} ) &\lesssim 2^{-(M+C)k'}
|E_3|^{(2-\theta)/2} |E_4|^{(\theta-1)/2}
\end{align*}
for some $0 < \theta < 1$ to be chosen later.

By Proposition \ref{abstract} we can thus bound the left-hand side of \eqref{weakin-2} by (if $M$ is chosen sufficiently large)
$$
2^{-k} 2^{-k'}
\frac{|E_1|^{(1+\theta_1)/2} |E_2|^{(1+\theta_2)/2}}{|E_4|^{1 - \theta_3}} (|E_3|^{(2-\theta)/2} |E_4|^{(\theta-1)/2})^{1-\theta_3},$$
and the claim follows by setting $\theta_1 := 2\alpha_1 - 1$, $\theta_2 := 2\alpha_2 - 1$, $\theta_3 := 2(\alpha_3 + \alpha_4) + 1$, and $\theta := (3\alpha_3 + 2\alpha_4) / (\alpha_3 + \alpha_4)$; the reader may verify that the constraints on $\theta_1, \theta_2, \theta_3, \theta$ can be obeyed for $\alpha$ arbitrarily close to $A_{1}, A_{2}$.

\section{Appendix: Proof of Proposition \ref{abstract}}\label{abstract-sec}

We now prove Proposition \ref{abstract}.  This is the analogue of Appendix III in \cite{mtt:walshbiest}, but using the necessary modifications for the Fourier case.  Our arguments shall be modeled on those in \cite{cct}, Section 9, which were in turn inspired by \cite{thiele'}.

Fix the collection $\Pv$ and the collections $a^{(j)}_{P_j}$.  We continue to assume
$\Pv$ is sparse.

We adopt the shorthand
$$ SIZE_j := \size_j( (a^{(j)}_{P_j})_{\pv \in \Pv} ); \quad 
ENERGY_j := \modenergy_j( (a^{(j)}_{P_j})_{\pv \in \Pv} ).$$

We may of course assume that $a^{(j)}_{P_j}$ are always non-zero.  We begin by considering the contribution of a single tree:

\begin{lemma}[Tree estimate]\label{single-tree}
Let $T$ be a tree in $\P$, and $a^{(j)}_{P_j}$ be complex numbers for all $P \in T$ and $j=1,2,3$.  Then
$$
| \sum_{\pv \in T} \frac{1}{|I_\pv|^{1/2}} a^{(1)}_{P_1} a^{(2)}_{P_2} a^{(3)}_{P_3}|
\leq |I_T|
\prod_{j=1}^3 \size_j( (a^{(j)}_{P_j})_{\pv \in T} ).$$
\end{lemma}

\begin{proof}  This proof is reproduced verbatim from \cite{mtt:walshbiest}.

Without loss of generality we may assume that $T$ is a 3-tree.  We then use H\"older to estimate the left-hand side by
$$ 
(\sum_{\pv \in T} |a^{(1)}_{P_1}|^2)^{1/2}
(\sum_{\pv \in T} |a^{(2)}_{P_2}|^2)^{1/2}
\sup_{\pv \in T} \frac{|a^{(3)}_{P_3}|}{|I_\pv|^{1/2}}.$$
From Definition \ref{size-def} we have
$$ (\sum_{\pv \in T} |a^{(j)}_{P_j}|^2)^{1/2} \leq |I_T|^{1/2} 
\size_j( (a^{(j)}_{P_j})_{\pv \in T} )$$
for $j=1,2$.  Also, since the singleton tree $\{\pv\}$ is a $1$-tree with top $\pv$, we have
$$ \frac{|a^{(3)}_{P_3}|}{|I_\pv|^{1/2}} \leq
\size_3( (a^{(3)}_{P_3})_{\pv \in T} )$$
for all $\pv \in \T$.  The claim follows.
\end{proof}

To bootstrap this summation over $T$ to a summation over $\Pv$ we would like to partition $\Pv$ into trees $T$ for which one has control over $\sum_T |I_T|$.  This will be accomplished by 

\begin{proposition}\label{decomp}  Let $1\le j\le 3$, 
$\Pv'$ be a subset of $\Pv$, $n \in \Z$, and suppose that
\be{size-stop}
\size_j( (a^{(j)}_{P_j})_{\pv \in \Pv'} ) \leq 2^{-n} ENERGY_j.
\end{equation}
Then we may decompose $\Pv' = \Pv'' \cup \Pv'''$ such that
\be{size-lower}
\size_j( (a^{(j)}_{P_j})_{\pv \in \Pv''} ) \leq 2^{-n-1}
ENERGY_j
\end{equation}
and that $\Pv'''$ can be written as the disjoint union of trees $\T$ such that
\be{tree-est}
\sum_{T \in \T} |I_T| \lesssim 2^{2n}.
\end{equation}
\end{proposition}

\begin{proof}
The idea is to initialize $\Pv''$ to equal $\Pv'$, and remove trees from $\Pv''$ one by one (placing them into $\Pv'''$) until \eqref{size-lower} is satisfied.

If $P$ is a tile, let $\xi_P$ denote the center of $\omega_P$.
If $P$ and $P'$ are tiles, we write $P' \lesssim^+ P$ if $P' \lesssim' P$ and
$\xi_{P'} > \xi_P$, and $P' \lesssim^- P$ if $P' \lesssim' P$ and $\xi_{P'} < \xi_P$.  

We now perform the following algorithm.  We shall need a collection $\T$ of trees, which we initialize to be the empty set.  We consider the set of all
trees $T$ of type $i\neq j$ in $\Pv$ which are ``upward trees'' in the 
sense that 
\be{plus}
P_j \lesssim^+ P_{T,j} \hbox{ for all } \pv \in T
\end{equation}
and which satisfy the size estimate
\be{t-plus}
\sum_{\pv \in T} |a^{(j)}_{P_j}|^2 
\geq 2^{-2n-3} |I_T|.
\end{equation}
If there are no trees obeying \eqref{plus} and \eqref{t-plus}, we terminate the algorithm.  Otherwise, we choose $T$ among all such trees so that 
the center $\xi_{T,j}$ of $\omega_{P_T,j}$ is maximal (primary goal), 
and that $T$ is maximal with respect to set inclusion (secondary goal).  Let $T'$ denote the $j$-tree
$$ T' := \{ \pv \in \Pv \backslash T: P_j \leq P_{T,j} \}.$$
We remove both $T$ and $T'$ from $\Pv$, and add them to $\T$.  Then one repeats the algorithm until we run out of trees obeying \eqref{plus} and \eqref{t-plus}.

Since $\Pv$ is finite, this algorithm terminates in a finite number of steps, producing trees $T_1, T'_1, T_2, T'_2, \ldots, T_M, T'_M$.  We claim that
the trees $T_1, \ldots, T_M$ produced in this manner are strongly $j$-disjoint.
It is clear from construction that $T_s \cap T_{s'} = \emptyset$ for all
$s \neq s'$; by the rank 1 assumption we thus see that
$P_j \neq P'_j$ for all $\pv \in T_s$, $\pv' \in T_{s'}$, $s \neq s'$.

Now suppose for contradiction that we had tri-tiles
$\pv \in T_s$, $\pv' \in T_{s'}$ such that 
$2\omega_{P_j} \subsetneq 2\omega_{P'_j}$
and $I_{P'_j} \subseteq I_{T_s}$.  From the sparseness assumption
we thus have $|\omega_{P'_j}| \geq 10^9 |\omega_{P_j}|$.  
Since $P_j \lesssim^+ P_{T_s,j}$ and 
$P'_j \lesssim^+ P_{T_{s'},j}$, we thus see that
$\xi_{P_{T_{s'}},j} < \xi_{P_{T_s},j}$.  
By our selection algorithm this implies that $s < s'$.

Also, since $|\omega_{P'_j}| \geq 10^9 |\omega_{P_j}|$, $I_{P'_j} \subseteq I_{T_s}$, and $P_j \lesssim P_{T_s,j}$ we see that $P'_j \leq P_{T_s,j}$.  
Since $s < s'$, this means that $\pv' \in T'_s$.  But
$T'_s$ and $T_{s'}$ are disjoint by construction, which is a contradiction.
Thus the trees $T_s$ are strongly $j$-disjoint.  From this, \eqref{t-plus}, \eqref{size-stop}, and Definition \ref{energy-def} we see that
$$ \sum_{s=1}^M |I_{T_s}| \lesssim 2^{2n}.$$
Since $T'_s$ has the same top as $T_s$, we may thus add all the $T_s$ and $T'_s$
to $\T$ while respecting \eqref{tree-est}.

Now consider the set $\Pv$ of remaining tri-tiles.  We note that
\be{plus-left}
\sum_{\pv \in T: P_j \lesssim^+ P_{T,j}} |a^{(j)}_{P_j}|^2
< 2^{-2n-3} |I_T|
\end{equation}
for all trees $T$ in $\Pv$, since otherwise the portion of $T$ which obeyed \eqref{plus} would be eligible for selection by the above algorithm.

We now repeat the previous algorithm, but replace $\lesssim^+$ by $\lesssim^-$
(so that the trees $T$ are ``downward-pointing'' instead of ``upward-pointing'') and select the trees $T$ so that the center $\xi_{T,j}$ is \emph{minimized} rather
than maximized.  This yields a further collection of  trees to add to $\T$ while still respecting \eqref{tree-est}, and the remaining collection of tiles $\Pv$ has the property that
\be{plus-right}
\sum_{\pv \in T: P_j \lesssim^- P_{T,j}} |a^{(j)}_{P_j}|^2 
< 2^{-2n-3} |I_T|
\end{equation}
for all trees $T$ in $\Pv$.  Combining \eqref{plus-left} and \eqref{plus-right} we obtain \eqref{size-lower} as desired.
\end{proof}

From Proposition \ref{decomp} we easily have

\begin{corollary}\label{decomp-cor}  There exists a partition
$$ \Pv = \bigcup_{n \in \Z} \Pv_n$$
where for each $n \in \Z$ and $j = 1,2,3$ we have
$$ \size_j( (a^{(j)}_{P_j})_{\pv \in \Pv'} ) \leq \min(2^{-n} ENERGY_j, SIZE_j).$$
Also, we may cover $\Pv_n$ by a collection $\T_n$ of trees such that
$$ \sum_{T \in \T_n} |I_T| \lesssim 2^{2n}.$$
\end{corollary}

\begin{proof}
Since $\Pv$ is finite, we see that the hypotheses of Proposition \ref{decomp} hold for all $1\le j\le 3$ if $n = -N_0$ for some sufficiently large $N_0$.  Set the $\Pv_n$ to be empty for all $n < -N_0$.  Now initialize $n = -N_0$ and $\Pv' = \Pv$.  For $1\le j\le 3$ in turn, we apply Proposition \ref{decomp}, moving the tri-tiles in $\Pv'''$ from $\Pv'$ in $\Pv_n$ and keeping the tri-tiles in $\Pv''$ inside $\Pv'$.  We then increment $n$ and repeat this process. Since we are assuming the $a^{(j)}_{P_j}$ are non-zero, every tri-tile must eventually be absorbed into one of the $\Pv_n$.  The properties are then easily verified.
\end{proof}

From Corollary \ref{decomp-cor} and Lemma \ref{single-tree} we see that
$$
| \sum_{\pv \in T} \frac{1}{|I_\pv|^{1/2}} a^{(1)}_{P_1} a^{(2)}_{P_2} a^{(3)}_{P_3}|
\leq |I_T|
\prod_{j=1}^3 \min(2^{-n} ENERGY_j, SIZE_j)$$
for all $T \in \T_n$.  Summing over all $T$ in $\T_n$ and then summing over all $n$, we obtain
$$
| \sum_{\pv \in T} \frac{1}{|I_\pv|^{1/2}} a^{(1)}_{P_1} a^{(2)}_{P_2} a^{(3)}_{P_3}|
\lesssim \sum_n 2^{2n}
\prod_{j=1}^3 \min(2^{-n} ENERGY_j, SIZE_j).$$
Without loss of generality we may assume that
$$ \frac{ENERGY_1}{SIZE_1} \leq \frac{ENERGY_2}{SIZE_2} \leq \frac{ENERGY_3}{SIZE_3}.$$
We may estimate the right-hand side as
\begin{align*}
 \sum_n \min(&2^n ENERGY_1 SIZE_2 SIZE_3, \\
& ENERGY_1 ENERGY_2 SIZE_3, \\
&2^{-n} ENERGY_1 ENERGY_2, ENERGY_3)
\end{align*}
which can be bounded by
$$ ENERGY_1 ENERGY_2 SIZE_3 \log(1 + \frac{ENERGY_3/SIZE_3}{ENERGY_2/SIZE_2}).$$
The claim then follows.

\end{document}